\documentclass[11pt,letterpaper]{article}

\usepackage[backend=bibtex]{biblatex} 
\usepackage{appendix}
\usepackage{longtable}
\usepackage{booktabs}
\usepackage{mathrsfs}

\usepackage{amsmath,amsthm,amsfonts,color,graphicx,xcolor,dsfont}
\usepackage{pdfsync,float,tcolorbox,mathabx,hyperref, halloweenmath}
\usepackage{amssymb}
\usepackage{array}
\usepackage{color}
\usepackage{multirow}
\usepackage{graphicx}
\usepackage{float} 
\usepackage{subfigure}

\usepackage{algorithm}  
\usepackage{algpseudocode}  
\usepackage{amsmath}  

\DeclareUnicodeCharacter{02B9}{'}

\numberwithin{equation}{section}
\hypersetup{
    colorlinks=true,
    linkcolor=blue,
    filecolor=magenta,      
    urlcolor=magenta,
  }

\begin{document}
	\title{Predicting band structures for 2D Photonic Crystals via Deep Learning}
\author{Yueqi Wang\thanks{Department of Mathematics, The University of Hong Kong, Pokfulam Road, Hong Kong. Email: {\tt{u3007895@connect.hku.hk}} YW acknowledges support from the Research Grants Council (RGC) of Hong Kong via the Hong Kong PhD Fellowship Scheme (HKPFS).}\and Richard Craster\thanks{Department of Mathematics, Imperial College London, London SW7 2AZ,UK. Email:\texttt{r.craster@imperial.ac.uk}} \and Guanglian Li\thanks{Corresponding author. Department of Mathematics, The University of Hong Kong, Pokfulam Road, Hong Kong. Email: {\tt{lotusli@maths.hku.hk}} GL acknowledges the support from Newton International Fellowships Alumni following-on funding awarded by The Royal Society, Young Scientists fund (Project number: 12101520) by NSFC and General Research Fund (GRF) (Project number: 17317022), RGC, Hong Kong.}}
	\maketitle
	\begin{abstract}
Photonic crystals (PhCs) are periodic dielectric structures that exhibit unique electromagnetic properties, such as the creation of band gaps where electromagnetic wave propagation is inhibited. Accurately predicting dispersion relations, which describe the frequency and direction of wave propagation, is vital for designing innovative photonic devices. However, traditional numerical methods, like the Finite Element Method (FEM), can encounter significant computational challenges due to the multiple scales present in photonic crystals, especially when calculating band structures across the entire Brillouin zone. To address this, we propose a supervised learning approach utilizing U-Net, along with transfer learning and Super-Resolution techniques, to forecast dispersion relations for 2D PhCs. Our model reduces computational expenses by producing high-resolution band structures from low-resolution data, eliminating the necessity for fine meshes throughout the Brillouin zone. The U-Net architecture enables the simultaneous prediction of multiple band functions, enhancing efficiency and accuracy compared to existing methods that handle each band function independently. Our findings demonstrate that the proposed model achieves high accuracy in predicting the initial band functions of 2D PhCs, while also significantly enhancing computational efficiency. This amalgamation of data-driven and traditional numerical techniques provides a robust framework for expediting the design and optimization of photonic crystals. The approach underscores the potential of integrating deep learning with established computational physics methods to tackle intricate multiscale problems, establishing a new benchmark for future PhC research and applications.
	\end{abstract}
\section{Introduction}
Photonic crystals (PhCs) are dielectric materials that are constructed from a unit cell periodically repeated with a period size comparable to the wavelength \cite{joannopoulos2008molding}. These materials exhibit the band gap phenomenon, where specific frequency intervals prevent electromagnetic wave propagation.
Complete band gaps occur when all polarizations and directions of wave propagation are prohibited, making PhCs promising for innovative photonic devices like optical transistors, photonic fibers, and low-loss optical mirrors 
 \cite{yanik2003all,russell2003photonic,wang2023analytical,labilloy1997demonstration}.
The band structure of PhCs describes how electromagnetic wave propagation depends on frequency, polarization, and direction, representing the dispersion relation mathematically.  In this paper, we focus on 2-dimensional (2D) periodic PhCs, which are homogeneous along the $z$ axis and have high-contrast dielectric materials embedded in dielectric materials within the $xy$ plane. We propose a novel supervised learning based scheme to achieve an accurate prediction of the dispersion relation for a given unit cell structure of 2D PhCs. 

Various numerical methods have been employed to calculate photonic band structures, including plane wave expansion methods \cite{ho1990existence,leung1990full}, the transfer matrix method \cite{pendry1992calculation,pendry1996calculating}, finite difference time domain method \cite{chan1995order}, layer Korringa-Kohn-Rostoker
method \cite{stefanou1992scattering} and multipole methods \cite{nicorovici1995photonic,botten2001photonic}. In particular, there is wide application of finite element methods in recent years \cite{axmann1999efficient,boffi2006modified,dobson1999efficient,engstrom2010complex,schmidt2009computation,schmidt2010efficient}. However, those numerical methods have to meet the issue of simulating the band structures with potentially high computational costs. The permittivity can vary widely, and the ratio between these values, the so-called contrast, should be large to generate the band gap. This leads to Helmholtz eigenvalue problems with both high-contrast and piecewise constant coefficients, which is numerically challenging. Thus, in the process of designing PhCs structures, it is desirable to establish a bidirectional relationship between the structure and the band gap property efficiently.

In the realm of computational mathematics, deep learning models offer a promising alternative to traditional approaches for modeling complex input-output mappings. These models have been successfully applied in various fields, including fluid dynamics \cite{beigzadeh2012prediction,butz2002modelling,mi2001flow} and materials science \cite{cang2018improving,koker2007neural,kondo2017microstructure}, showcasing their potential in capturing elusive relationships between material structures and properties. Unsurprisingly, deep learning models have been already utilized in the dispersion relation prediction and the inverse design of PhCs \cite{liu2018training,jiang2022dispersion,yao2019intelligent,tahersima2019deep,10.1093/jcde/qwad013}. 
In particular, Jiang et al. \cite{jiang2022dispersion} propose the use of convolutional neural networks (CNNs) and conditional generative adversarial networks (cGANs) to bridge the structure and properties from the
forward and inverse directions, respectively. Recent work by Ma et al. \cite{PhysRevX.11.021052} also demonstrates the potential of neural networks in predicting photonic band structures, showcasing how data-driven methods can achieve remarkable efficiency in the exploration of PhC designs. 


Unlike those existing approaches, we propose a novel supervised learning scheme which leverages the U-Net architecture with transfer learning and Super-Resolution techniques to predict dispersion relations efficiently. This integration of data-driven and traditional computational methods holds the potential to accelerate progress in designing and optimizing PhCs structures. The main features of our scheme are fourfold. 
(a) It is capable of predicting the entire dispersion relation, or the first several band functions, using one supervised learning scheme. Current deep learning based schemes either use separate neural networks to predict each band function or use one neural network with all the band functions of interest as output \cite{jiang2022dispersion,Christensen2020}. (b) It predicts dispersion relation over the entire first Brillouin zone instead of only on its boundary, and hence is more accurate.
(c) By hybridizing the U-Net architecture with transfer learning, our supervised learning scheme can learn effectively from relatively small datasets and reduce the overall training cost significantly \cite{chollet2021deep}.
(d) By incorporating a Super-Resolution (SR) model, our scheme has the potential to reduce the cost of generating training data or improve the accuracy of existing noisy data. This represents a huge reduction of computational complexity since traditional numerical methods suffer from high computational cost, which makes the generation of the training data extremely expensive. 

The remainder of the paper is organized as follows. In Section \ref{Problem formulation}, we revisit the mathematical models used for calculating band functions and outline the primary supervised learning tasks. Then we delve into the neural network architectures employed in our approach in Section \ref{sec:neuralnetwork}. The training and testing procedures as well as the performance and accuracy of the proposed model are demonstrated in Section \ref{Results}. Finally, we consolidate our findings and conclusions in Section \ref{sec:conclusion}.

 
\section{Problem formulation}\label{Problem formulation}
In this section, we provide a recap of the mathematical framework used for calculating band functions and introduce the primary supervised learning tasks that form the foundation of our approach.
\subsection{Dispersion relation}

In the SI convention, the time harmonic Maxwell equations for linear, non-dispersive, and nonmagnetic media, with free charges and free currents, consist of a system of four equations  \cite{jackson1999classical},
\begin{subequations}
	\begin{align}
		\nabla\times\mathbf{E}(\mathbf{x})-i\omega\mu_{0}\mathbf{H}(\mathbf{x})=0,\label{harm1}    \\
		\nabla\times\mathbf{H}(\mathbf{x})+i\omega\epsilon_{0}\epsilon(\mathbf{x})\mathbf{E}(\mathbf{x})=0 , \label{harm2}\\  \nabla\cdot\left(\epsilon(\mathbf{x})\mathbf{E}(\mathbf{x})\right)=0,\label{harm3}\\
		\nabla\cdot\mathbf{H}(\mathbf{x})=0.\label{harm4}
	\end{align}
\end{subequations}
Here $\mathbf{x}\in\mathbb{R}^3$, $\mathbf{E}$ is the electric field, $\mathbf{H}$ the magnetic field and $\mathbf{D}$ the electric displacement field. The scalar $\omega\geq 0$ is the frequency of the electromagnetic wave, $\mu_0$ the vacuum permeability, $\epsilon_0$ the vacuum permittivity, and $\epsilon\in L^{\infty}(\mathbb{R}^3;\mathbb{R}^+)$ the relative permittivity.

Applying the curl operator to \eqref{harm1} and using \eqref{harm2}, we obtain
\begin{equation}\label{E}
	\nabla\times\left(\nabla\times\mathbf{E}(\mathbf{x})\right)-(\omega c^{-1})^2\epsilon(\mathbf{x})\mathbf{E}(\mathbf{x})=0,
\end{equation}
with $\epsilon_{0}\mu_{0}=c^{-2}$, where $c$ is the speed of light. Similarly, applying 
the curl operator to \eqref{harm2} and using \eqref{harm1}, we obtain
\begin{equation}\label{H}
	\nabla\times\left(\epsilon(\mathbf{x})^{-1}\nabla\times\mathbf{H}(\mathbf{x})\right)-\left(\omega c^{-1}\right)^2\mathbf{H}(\mathbf{x})=0.
\end{equation}

In this study, our focus is on 2D PhCs, which exhibit finite extensions in the $z$ direction and finite periodicities in the $x$-$y$ plane. However, it is typically assumed that the material extends infinitely in the plane perpendicular to the columns. Under such setting, the relative permittivity $\epsilon(\mathbf{x})$ is assumed to be independent of the $z$ direction. Then we can split the electromagnetic fields $\mathbf{E}=(E_1,E_2,E_3)$ and $\mathbf{H}=(H_1,H_2,H_3)$ in \eqref{E} and \eqref{H} into the transverse electric (TE) mode with $H_1=H_2=E_3=0$ and transverse magnetic (TM) mode with $E_1=E_2=H_3=0$. Each mode is a scalar eigenvalue problem,
\begin{align}
	-\nabla \cdot (\epsilon(\mathbf{x})^{-1}\nabla H(\mathbf{x}))-(\omega c^{-1})^2 H(\mathbf{x})&=0,\quad \mathbf{x}
	\in\mathbb{R}^2,&&\textbf{(TE mode)},\label{TE}\\
	-\Delta E(\mathbf{x})-(\omega c^{-1})^2 \epsilon(\mathbf{x})E(\mathbf{x})&=0,\quad \mathbf{x}
	\in\mathbb{R}^2.&& \textbf{(TM mode)}.\label{TM}
\end{align}
Since 2D PhCs possess a discrete translational symmetry in the $x$-$y$ plane \cite{joannopoulos2008molding}, the relative
permittivity $\epsilon(\mathbf{x})$ satisfies
\begin{equation*}
	\epsilon(\mathbf{x}+c_1\mathbf{a}_1+c_2\mathbf{a}_2)=\epsilon(\mathbf{x}),\quad
	\forall\mathbf{x}\in\mathbb{R}^2\text{ and }c_1,c_2\in\mathbb{Z}.
\end{equation*}
The primitive lattice vectors, denoted by $\mathbf{a}_1$ and $\mathbf{a}_2$, are the shortest possible vectors that fulfill this condition and they span the fundamental periodicity domain $\Omega$, also known as unit cell. The reciprocal lattice vectors are defined by
\begin{equation}
	\mathbf{b}_i \cdot \mathbf{a}_j=2\pi\delta_{ij},\quad \text{ for }i,j=1,2,
\end{equation}
which generate the so-called reciprocal lattice. The elementary cell of the reciprocal lattice is the (first) Brillouin zone $\mathcal{B}_F$, i.e., the region closer to a certain lattice point than to any other lattice points in the reciprocal lattice. In particular, in this work we consider a square lattice, whose primitive lattice vectors are $\mathbf{a}_i = a\mathbf{e}_i$ for $i=1,2$. Here, $(\mathbf{e}_i)_{i=1,2}$ is the canonical basis in $\mathbb{R}^2$ and $a\in \mathbb{R}^+$ lattice constant. The corresponding reciprocal lattice vectors are $\mathbf{b}_i = \frac{2\pi}{a}\mathbf{e}_i$ for $i=1,2$.

\begin{figure}[hbt!]
	\centering
	\subfigure[Square lattice]{\label{lattice(1)}
		\includegraphics[width = .27\textwidth,trim={3cm 1cm 2.5cm 0.8cm},clip]{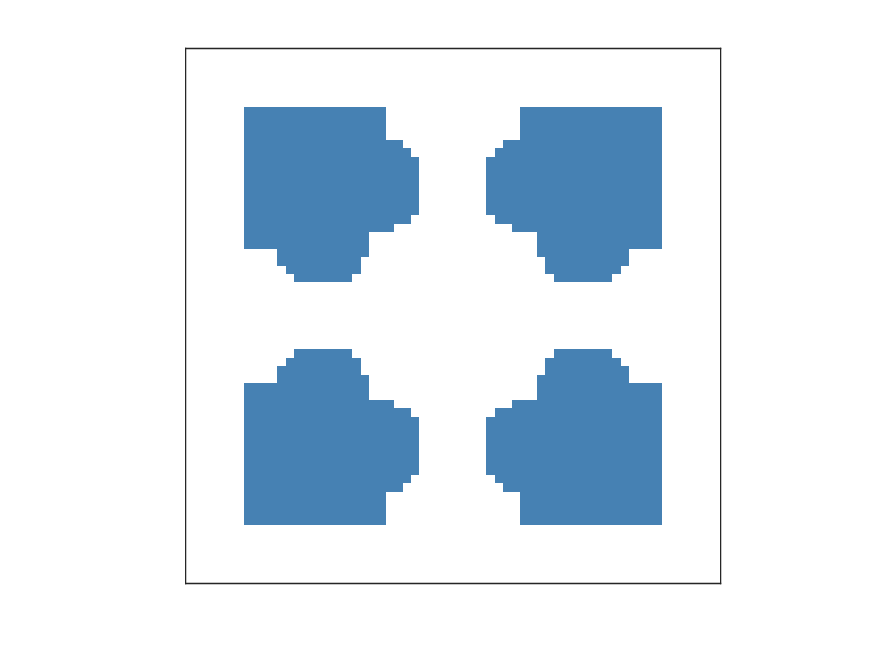}
	}%
	\subfigure[The corresponding Brillouin zone]{\label{lattice(2)}
		\includegraphics[width = .33\textwidth,trim={16cm 5.8cm 14cm 4.5cm},clip]{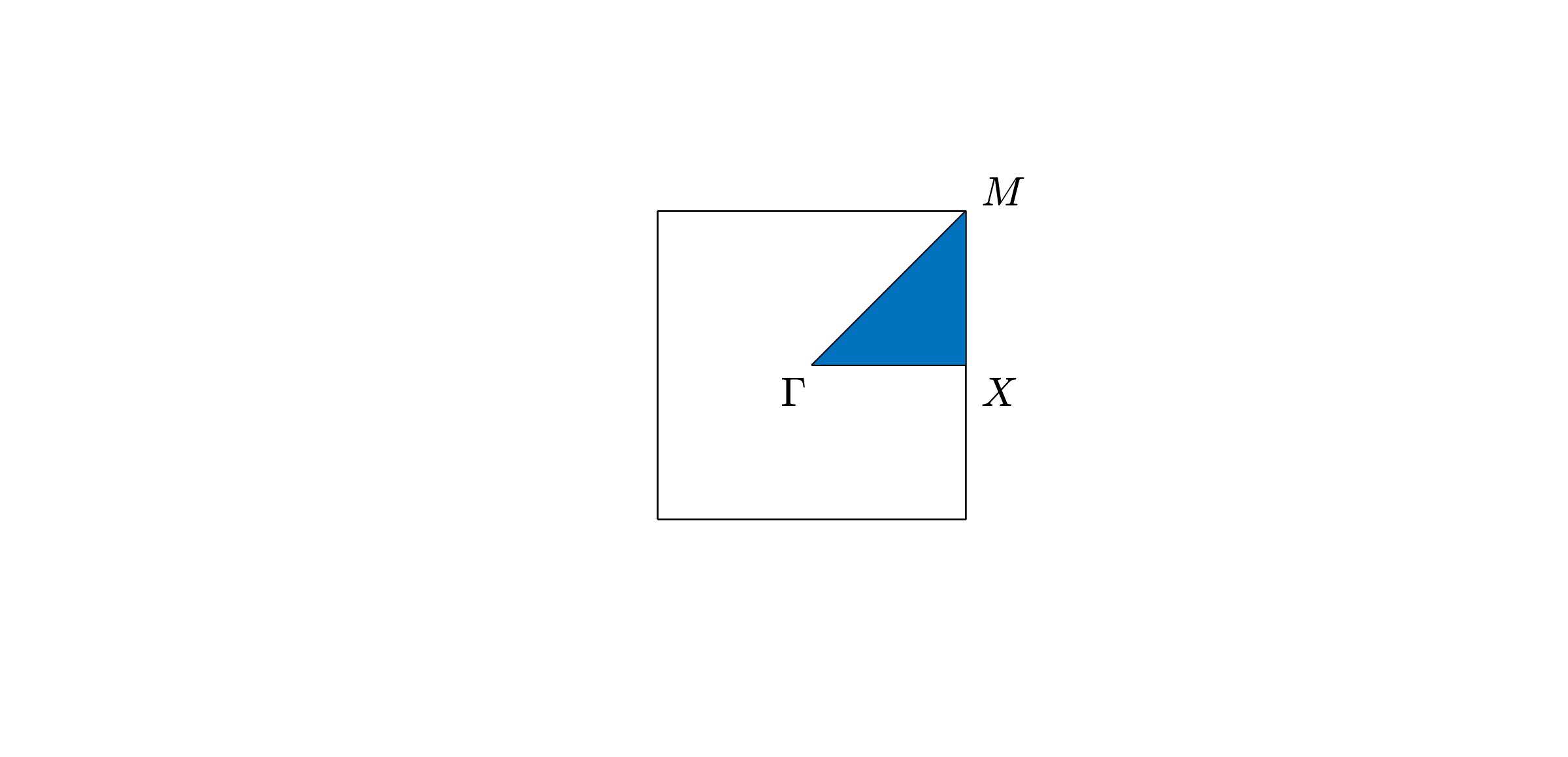}
	}%
	\centering
	\caption{Illustration of unit cell $\Omega$ for a square lattice (left) and the corresponding Brillouin zone (right): In $\Omega$, blue stands for alumina with permittivity $8.9$ and white for air with permittivity $1$; In $\mathcal{B}_F$, the IBZ is the shaded triangle with vertices $\Gamma=(0,0)$, $X=\frac{1}{a}(\pi,0)$ and $M=\frac{1}{a}(\pi,\pi)$.}\label{lattice}
\end{figure}

Bloch's theorem \cite{kittel2018introduction} states that in periodic crystals, wave functions take the form of a plane wave modulated by a periodic function. Thus they can be written as
$\Psi(\mathbf{x})=e^{i\mathbf{k}\cdot\mathbf{x}}u(\mathbf{x})$,
where $\Psi$ is the wave function, $u(\mathbf{x})$ is periodic sharing the periodicity of the crystal lattice and $\mathbf{k}$ is the wave vector varying in the  Brillouin zone $\mathcal{B}_F$. The periodic condition of $u(\mathbf{x})$ implies that each wave function is determined by its values within the unit cell $\Omega$. Thus, the solutions to \eqref{TE} and \eqref{TM} are expressed as
$H(\mathbf{x})=e^{i\mathbf{k}\cdot \mathbf{x}}u_1(\mathbf{x})$ and $E(\mathbf{x}) =e^{i\mathbf{k}\cdot \mathbf{x}}u_2(\mathbf{x})$ for some periodic functions $u_1(\mathbf{x})$ and $u_2(\mathbf{x})$ in the unit cell $\Omega$, and the parameterized Helmholtz eigenvalue problems \eqref{TE} and \eqref{TM} reduce to
\begin{subequations}
	\begin{align}
		-(\nabla+i\mathbf{k})\cdot\left(\epsilon(\mathbf{x})^{-1}(\nabla+i\mathbf{k}) u_1(\mathbf{x})\right)-\left( \omega c^{-1}\right)^{2} u_1(\mathbf{x})&=0,\quad \mathbf{x}\in\Omega && \textbf{(TE mode)},\label{TE2}\\
		-(\nabla+i\mathbf{k})\cdot\left((\nabla+i\mathbf{k}) u_2(\mathbf{x})\right)-\left(\omega c^{-1}\right)^{2}\epsilon(\mathbf{x})u_2(\mathbf{x})&=0,\quad \mathbf{x}\in\Omega&& \textbf{(TM mode)},\label{TM2}
	\end{align}
\end{subequations}
where $\mathbf{x}\in\Omega\subset\mathbb{R}^2$, $\mathbf{k}$ varies in the Brillouin zone (BZ), and $u_i(\mathbf{x})$ satisfies
the periodic boundary conditions $u_i(\mathbf{x})=u_i(\mathbf{x}+\mathbf{a_j})$ with $\mathbf{a_j}$ being the primitive lattice vector for $i,j=1,2$.
If the materials in the unit cell have additional symmetry, e.g. mirror symmetry, we can further restrict $\mathbf{k}$ to the triangular irreducible Brillouin zone (IBZ), denoted by $\mathcal{B}$. An example of the square lattice and its Brillouin zone are shown in Figures \ref{lattice(1)} and \ref{lattice(2)}.

In summary, we can formulate both parameterized Helmholtz problems \eqref{TE2} and \eqref{TM2} by
\begin{equation}\label{both}
	-(\nabla+i\mathbf{k})\cdot \alpha(\mathbf{x})(\nabla+i\mathbf{k})u(\mathbf{x})-\lambda\beta(\mathbf{x})u(\mathbf{x})=0, \quad \mathbf{x}\in\Omega,
\end{equation}
with $\Omega\subset\mathbb{R}^2$, $\mathbf{k} \in \mathcal{B}$, and $\lambda =(\omega c^{-1})^2$. In the TE mode, $u$ describes the magnetic field $H$ in the $z$-direction and the coefficients $\alpha(\mathbf{x})$ and $\beta(\mathbf{x})$ are $ \alpha(\mathbf{x}):=\epsilon(\mathbf{x})^{-1}$ and $\beta(\mathbf{x}):=1$.
Similarly, in the TM mode, $u$ describes the electric field $E$ in the $z$-direction and the coefficients $\alpha(\mathbf{x})$ and $\beta(\mathbf{x})$ are
$ \alpha(\mathbf{x}):=1$ and $\beta(\mathbf{x}):=\epsilon(\mathbf{x})$.

The variational formulation of \eqref{both} reads: find a non-trivial eigenpair $(\lambda,u) \in (\mathbb{R},H^1_{\pi}(\Omega))$ for $\mathbf{k}\in\mathcal{B}$ such that
\begin{equation}\label{variational}
	\left\{
	\begin{aligned}   \int_{\Omega}\alpha(\nabla+i\mathbf{k})u\cdot(\nabla-i\mathbf{k})\Bar{v}-\lambda\beta u\Bar{v}\mathrm{d}x&=0, \quad \forall v \in H^1_{\pi}(\Omega)\\
		\|u\|_{L^2_{\beta}(\Omega)}&=1.
	\end{aligned}
	\right.
\end{equation}
Here, we define $L^2_{\beta}(\Omega)$ as the space of weighted square integrable
functions equipped with the norm
$\|f\|_{L^2_{\beta}(\Omega)}:=(\int_{\Omega} \beta(\mathbf{x})|f(\mathbf{x})|^2\mathrm{d}\mathbf{x})^{\frac{1}{2}}$.
Let $H^1(\Omega)\subset L^2_{\beta}(\Omega)$ be with square integrable gradient, equipped with the standard $H^1(\Omega)$-norm and $H^1_{\pi}(\Omega)\subset H^1(\Omega)$ is composed of functions with periodic boundary conditions.
Using the sesquilinear forms
\begin{align*}
	a(u,v)&:=\int_{\Omega}\alpha (\nabla+i\mathbf{k})u\cdot(\nabla-i\mathbf{k})\Bar{v}\,\mathrm{d}x\quad \mbox{and}\quad
	b(u,v):=\int_{\Omega}\beta u\Bar{v}\,\mathrm{d}x,
\end{align*}
problem \eqref{variational} is equivalent to finding a non-trivial eigenpair $(\lambda,u) \in (\mathbb{R},H^1_{\pi}(\Omega))$ for $\mathbf{k}\in\mathcal{B}$ such that
\begin{equation}\label{simply}
	\left\{
	\begin{aligned}
		a(u,v)&=\lambda b(u,v), \text{ for all }v \in H^1_{\pi}(\Omega) \\
		\|u\|_{L^2_{\beta}(\Omega)}&=1.
	\end{aligned}
	\right.
\end{equation}

This eigenvalue problem \eqref{simply} has a non-negative discrete eigenvalue sequence which can be enumerated in a nondecreasing manner and repeat according to their finite multiplicities as \cite{glazman1965direct}
\begin{align*}  0\leq\lambda_1(\mathbf{k})\leq\lambda_2(\mathbf{k})\leq\cdots \leq \lambda_n(\mathbf{k})\leq\cdots\leq \infty.
\end{align*}
$\{\lambda_n(\mathbf{k})\}_{n=1}^{\infty}$ is an infinite sequence with $\lambda_n(\mathbf{k})$ being a continuous function with respect to the wave vector $\mathbf{k}$, $\lambda_n(\mathbf{k})\to \infty$ when $n \to \infty$ and the graph of $\{\lambda_n(\mathbf{k})\}_{n=1}^{\infty}$ is a real analytic variety in $\mathbb{R}^3$ \cite{kuchment1993floquet}. 

Subsequently, the $n$th band function $\omega_n(\mathbf{k}) := c\sqrt{\lambda_n(\mathbf{k})}$ is considered a function of the wave vector $\mathbf{k}$ for all $n\in\mathbb{N}^+$. This representation allows for the visualization of a graph where each band function $\omega_n$ is plotted against the wave vector $\mathbf{k}$. The regions between adjacent band functions are known as band gaps, and the dispersion relation describes the relationship between $\omega_n$ and $\mathbf{k}$.

\subsection{Supervised learning tasks}
We propose a supervised learning-based method to predict the first $L$ band functions, where $L\in\mathbb{N}^+$. As the graph of ${\lambda_n(\mathbf{k})}_{n=1}^{\infty}$ forms a real analytic variety, there exists potential correlation among different band functions. This motivates us to use a single network as the approximator for the first $L$ band functions simultaneously. Therefore, the input of our neural networks contains information on the cell of the photonic crystal, the band number $n$, and the output provides an approximation to $\omega_n$ for $n=1,\cdots,L$.

\begin{figure}[htp!]
	\centering
	\includegraphics[width = .75\textwidth,trim={2cm 2.5cm 2cm 1.5cm},clip]{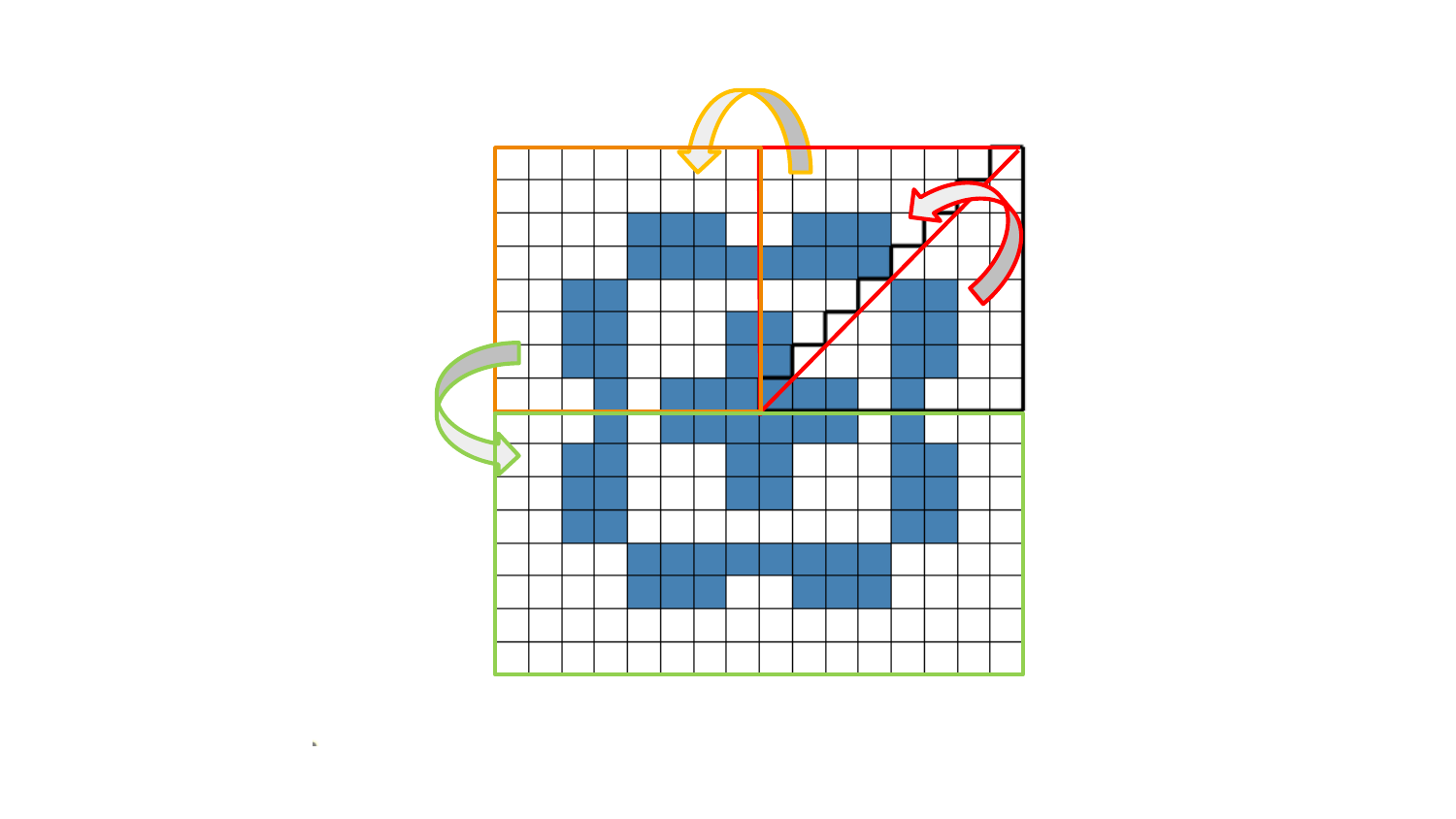}
	\centering
	\caption{Example of unit cell with $p4m$ plane symmetry}\label{generate_UnitCell}
\end{figure}

Our focus is on 2D PhCs with square unit cells that exhibit $p4m$ plane symmetry, as illustrated in Figure \ref{generate_UnitCell}. We specifically consider the transverse electric (TE) mode for 2D PhCs, described by equation \eqref{H}. Following the approach in \cite{jiang2022dispersion}, our dataset comprises $\tilde{N}=100,000$ unit cells with $p4m$ symmetry at resolutions of $m \times m$, where $m = 64, 16$.

For any given unit cell $\Omega$, we first discretize it by a binary image with a given size of $m\times m$, denoted as $M_{\Omega,m}$, such that
\begin{equation}\label{eq:dis-ucell}
	M_{\Omega,m}(i,j)=\left\{
	\begin{aligned}
		&0\\
		&1.
	\end{aligned}
	\right.
\end{equation}
Here, 0 stands for alumina with permittivity $8.9$ and 1 for air with permittivity $1$. 
For any given band number $L\in\mathbb{N}$, we generate a constant matrix of the same size $I_m^{(n)}$ with $I_m^{(n)}(i,j)=n/L$ for all $i,j=1,\cdots,m$, $n=1,\cdots,L$. 
Next, we consider two different tasks in the supervised learning. The first task is to learn a mapping from the discrete unit cell structure to any band function $F: {\{0,1\}}^{m^2}\times \mathbb{R}^{m^2}\to \mathbb{R}^{m^2}$,
\begin{align}\label{F 1}
	F:M_{\Omega,m}\times I_m^{(n)}\mapsto \hat{\omega}_n^{m}, 
\end{align}
where the output matrix $\hat{\omega}_n^{m}$ in \eqref{F 1} represents the prediction of the $n$th band function $\omega_n^{m}$ for a given unit cell $\Omega$, calculated over $m \times m$ uniformly distributed points in the first Brillouin zone.

In the second task, we aim to predict band functions with high accuracy from approximate band functions with low accuracy. This input can arise from the outcome of learning task \eqref{F 1} with smaller $m$. Hence, we want to learn a mapping from "low-resolution" band function to "high-resolution" one $G: \mathbb{R}^{m^2_1}\to \mathbb{R}^{m^2_2}$ with $m_1<m_2$,
\begin{align}\label{F 2}
		G:\omega_n^{m_1}\mapsto \hat{\omega}_n^{m_2}.
\end{align}
Here,  $\hat{\omega}_n^{m_2}$ in \eqref{F 2} is the predicted result of enhancing the resolution of $\omega_n^{m_1}$ from $m_1 \times m_1$ uniformly distributed points to $m_2 \times m_2$.


\section{Neural network}\label{sec:neuralnetwork}
In this section, we introduce the neural network architectures for the supervised learning tasks \eqref{F 1} and \eqref{F 2}.
\subsection{U-Net Network with transfer learning}\label{UNET}

We propose to use a U-Net architecture \cite{ronneberger2015u}, as the surrogate model for $F$ \eqref{F 1}, which has been extensively used for image-to-image translation tasks \cite{farimani2017deep,thuerey2020deep,wang2022deep,fotiadis2020comparing,isola2017image} and has shown promising performance in efficiently handling the mapping between high-dimensional inputs and outputs. 

The U-Net architecture we employed consists of a series of two consecutive convolutional operations, followed by either max-pooling or transposed convolution operations in the encoder and decoder pathways, respectively. In the encoder, we use four convolutional blocks to transform an initial $m\times m\times 2$ input matrix $M_{\Omega,m}\times I_m^{(n)}$, denoted by $\mathbf{F}_1^{(0)}\in\mathbb{R}^{m\times m\times2}$, into a  feature map $\mathbf{F}_{4}^{(2)}\in \mathbb{R}^{\frac{m}{16}\times \frac{m}{16}\times 1024}$. This feature map is then processed through two additional 2D convolutions, extending the output to 2048 feature channels, denoted by $\mathbf{Z}_{5}^{(2)}$, which is subsequently fed into the decoder. 

\begin{figure}[htp!]
	\centering
	\includegraphics[width = 0.9\textwidth,trim={0cm 1cm 0cm 1cm},clip]{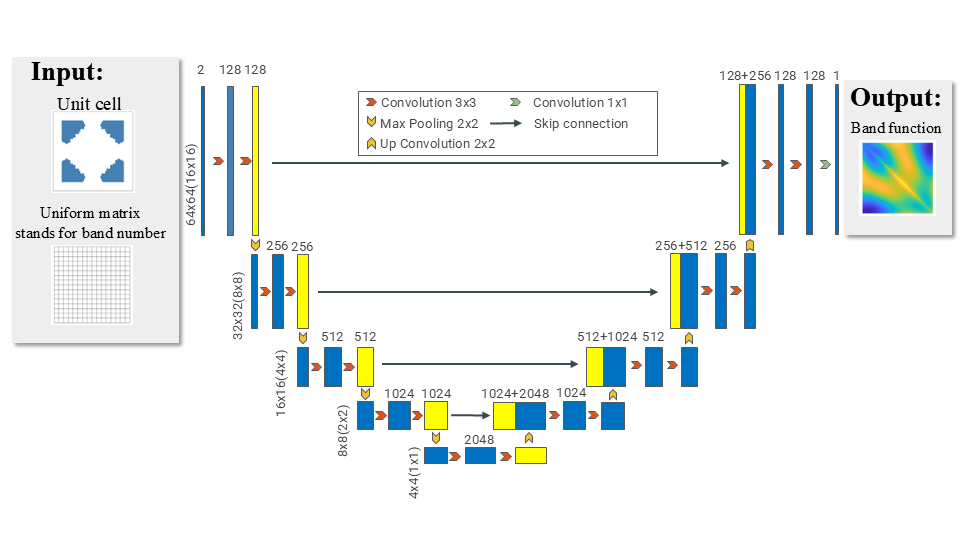}
	\centering
	\caption{The workflow of predicting dispersion relation with given unit cell.}\label{unet structure}
\end{figure}

Each convolutional block in the encoder consists of two consecutive convolutional operations, represented as follows:
\begin{align*}
	\mathbf{F}^{(i+1)}_{\ell}=\sigma\left({\rm BN}\left(\mathbf{W}_{\ell,\text{en}}^{(i)}*\mathbf{F}_{\ell}^{(i)}+\mathbf{b}_{\ell,\text{en}}^{(i)}\right)\right),\quad i=0,1,\quad \ell=1,\cdots,4,
\end{align*}
where $*$ denotes convolution, $\sigma$ represents the ReLU activation function, ${\rm BN}$ is short for the batch normalization, $\mathbf{W}_{\ell,\text{en}}^{(i)}$ is the $i$-th convolutional kernel for layer $\ell$, and $\mathbf{b}_{\ell,\text{en}}^{(i)}$ denotes the corresponding bias. Here, each $\mathbf{W}_{\ell,\text{en}}^{(i)}$ is of dimension $3\times 3\times d_{\ell,\text{en},\text{in}}^{(i)}\times d_{\ell,\text{en},\text{out}}^{(i)}$, where $d_{\ell,\text{en},\text{in}}^{(i)}$ and $d_{\ell,\text{en},\text{out}}^{(i)}$ are the input and output feature channels, respectively. We consider $d_{1,\text{en},\text{in}}^{(0)}=2$ and $d_{\ell,\text{en},\text{in}}^{(0)}=2^{5+\ell}$ for $\ell=2,3,4$, $d_{\ell,\text{en},\text{out}}^{(0)}=d_{\ell,\text{en},\text{in}}^{(1)}= d_{\ell,\text{en},\text{out}}^{(1)}=2^{6+\ell}$ for $\ell=1,\cdots,4$ in our numerical test as shown in Figure \ref{unet structure}. Note that here and in the following, we omit padding in the convolution operations within the formulas. Zero-padding is applied to ensure that the output dimensions meet our intended size. After the two consecutive convolutions, downsampling is performed using max-pooling, reducing the spatial dimensions by a factor of two. For a feature map $\mathbf{F}^{(2)}_{\ell} \in \mathbb{R}^{\frac{m}{2^{\ell-1}}\times\frac{m}{2^{\ell-1}} \times d}$ at layer $\ell$ in the encoder, max-pooling is defined as:
\begin{align*}
	\mathbf{F}^{(2)}_{\ell}(x,y,d)=\max_{(u,v)\in\{0,1\}^2}\mathbf{F}^{(2)}_{\ell}(2x+u,2y+v,d),\quad \ell=1,\cdots,4.
\end{align*}
Then, $\mathbf{F}^{(2)}_{\ell} \in \mathbb{R}^{\frac{m}{2^{\ell}} \times \frac{m}{2^{\ell}} \times d}$. This operation preserves the depth while halving both the height and width of the feature map.
In the decoder, symmetric upsampling is applied in four stages, eventually reconstructing a matrix $\mathbf{Z}_{1}^{(2)}\in \mathbb{R}^{m\times m\times128}$. The final output $\hat{\omega}_n^m:=\mathbf{Z}_{output}$ is then obtained by one $1\times 1$ convolution to reduce the number of feature channels back to the desired output dimension. The upsampling process at each stage is defined as
\begin{align*}
	\mathbf{Y}_{\ell}={\rm ConvTranspose}\left(\mathbf{Z}_{\ell+1}^{(2)},\mathbf{W}_{\ell,\text{up}},\mathbf{b}_{\ell,\text{up}}\right),\quad \ell=4,\cdots,1,
\end{align*}
where ${\rm ConvTranspose}$ denotes the transposed convolution operation with kernel $\mathbf{W}_{\ell,\text{up}}$ of size $2\times 2\times \left(2d_{\ell,\text{de},\text{out}}^{(0)}\right)\times \left(2d_{\ell,\text{de},\text{out}}^{(0)}\right)$ and bias $\mathbf{b}_{\ell,\text{up}}$. 

Next, skip connections are employed to concatenate feature maps from corresponding encoder and decoder layers, preserving low-level information for better reconstruction:
\begin{align*}
	\mathbf{Z}^{(0)}_{\ell}=\left[\mathbf{F}^{(2)}_{\ell},\mathbf{Y}_{\ell}\right]\in\mathbb{R}^{\frac{m}{2^{\ell-1}}\times \frac{m}{2^{\ell-1}}\times\left(3d_{\ell,\text{en},\text{out}}^{(1)}\right)},\quad \ell=4,\cdots,1,
\end{align*}
where $\mathbf{F}^{(2)}_{\ell}$ and $\mathbf{Y}_{\ell}$ are feature maps from the encoder and decoder at level $\ell$. This concatenation enhances feature extraction during reconstruction by retaining crucial spatial information.

Following this, two additional convolutional operations are applied: 
\begin{align*}
	\mathbf{Z}^{(i+1)}_{\ell}=\tilde{\sigma}\left(\mathbf{W}_{\ell,\text{de}}^{(i)}*\mathbf{Z}_{\ell}^{(i)}+\mathbf{b}_{\ell,\text{de}}^{(i)}\right),\quad i=0,1,\quad \ell=4,\cdots,1.
\end{align*}
Here, each $\mathbf{W}_{\ell,\text{en}}^{(i)}$ is of size $3\times 3\times d_{\ell,\text{de},\text{in}}^{(i)}\times d_{\ell,\text{de},\text{out}}^{(i)}$, and $d_{\ell,\text{de},\text{in}}^{(0)}=2^{6+\ell}+2^{7+\ell}$, $d_{\ell,\text{de},\text{out}}^{(0)}=d_{\ell,\text{de},\text{in}}^{(1)}= d_{\ell,\text{de},\text{out}}^{(1)}=2^{6+\ell}$ for $\ell=1,\cdots,4$ in our numerical test as shown in Figure \ref{unet structure}. 
This process doubles the spatial dimensions, reversing the downsampling effect from the encoder.

Considering the complexity of the U-Net structure, reducing computational load and speeding up training are crucial objectives. To address this, we incorporate transfer learning \cite{chollet2021deep} with the U-Net. Transfer learning involves using a pre-trained model as a starting point to accelerate learning on a new, but related task. In our numerical test where $L=10$, we initially train the model on band functions for $n=1,\cdots,5$ to establish a baseline model. Subsequently, we fine-tune the model on band functions for $n=6,\cdots,10$ using a smaller set of training samples. The two-stage training process we employ is efficient and capitalizes on the advantages of transfer learning, enhancing the robustness and computational efficiency of our approach. In contrast to the approach in \cite{jiang2022dispersion}, where individual networks are introduced for each band function, our model can predict multiple band functions concurrently, eliminating the necessity to train multiple independent models.
The workflow of the U-Net architecture with transfer learning for predicting dispersion relations is illustrated in Figure \ref{unet structure}.

Let $\{(M_{\Omega_j,m}, I_m^{(n)}),\omega_n^{m,j})\}_{j=1,\cdots,\tilde{N}_{\rm Train},n=1,\cdots,L}$ denote the labeled dataset, $N_{\rm Train}:=L\tilde{N}_{\rm Train}$ the total number of training data samples and $\theta_1$ the trainable parameters in the U-Net. To train the U-net \eqref{F 1}, we utilize the loss function based on the discrete mean square error 
\begin{align}\label{eq:loss-f1}
	\mathcal{L}_1(\theta_1)=\frac{1}{m^2N_{\rm Train}}\sum_{j=1}^{\tilde{N}_{\rm Train}}\sum_{n=1}^L\left\|F(M_{\Omega_j,m}, I_m^{(n)};\theta_1)
	-\omega_n^{m,j}\right\|_F^2.
\end{align}
Here, $\|\cdot\|_F$ denotes the Frobenius norm. The dataset is shuffled and the training is performed using batches of size $N_B$. 



\subsection{Super-Resolution Residual Network}\label{Super-Resolution Residual Network}

To enhance the resolution of the dispersion relation, we employ a method similar to image super-resolution. Unlike conventional upsampling techniques such as nearest neighbor, bilinear, or bicubic interpolation, our approach focuses on recovering finer details that may be lost when the dispersion relation is calculated on a coarse mesh in the Brillouin zone. We treat the band function values at $m_2\times m_2$ sampling points as "high-resolution" images and those at $m_1\times m_1$ sampling points as "low-resolution" images with $m_2:=4\times m_1$ and $m_1\in\mathbb{N}$. This dataset is used to train a neural network that can effectively increase the resolution of the dispersion relation.
\begin{figure}[htp!]
	\centering
	\includegraphics[width = .85\textwidth,trim={0cm 0cm 0cm 0cm},clip]{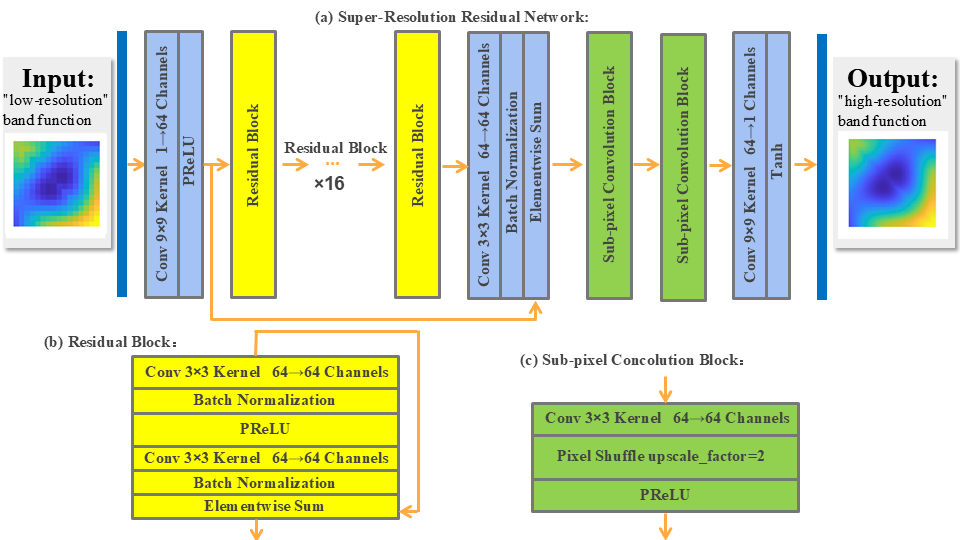}
	\centering
	\caption{The workflow of predicting dispersion relation with given "low-resolution" band function.}\label{SRResNet structure}
\end{figure}

For this task, we utilize the Super-Resolution Residual Network (SRResNet) \cite{ledig2017photo}, a generative fully convolutional neural network designed to enhance image resolution. The SRResNet magnifies the image resolution by a factor of 4 in each dimension, resulting in a 16-fold increase in the number of pixels. Comprising a series of residual blocks and sub-pixel convolution blocks, the SRResNet facilitates superior-quality super-resolution outcomes.

The initial "low-resolution" input matrix, $\mathbf{X}_{\text{low}}\in\mathbb{R}^{m_1\times m_1}$, is passed through a convolutional layer with kernel size $9 \times 9$ and stride $1$ to produce a feature map $\mathbf{F}_1^{(0)}\in\mathbb{R}^{m_1\times m_1 \times m_2}$:
\begin{align*}
\mathbf{F}_1^{(0)}={\rm PReLU}\left(\mathbf{W}_0*\mathbf{X}_{\text{low}}+\mathbf{b}_0\right),
\end{align*}
where $\mathbf{W}_{0}$ is of dimension $9\times 9 \times 1 \times m_2$, $\mathbf{b}_{0}$ is a bias term and a parametric ReLU (PReLU) activation is applied here. Next, $\mathbf{F}_1^{(0)}$ is passed through a sequence of 16 residual blocks while the resolution is maintained. Each residual block consists of two convolutional layers with kernel size $3\times3$ and stride $1$, followed by Batch Normalization (BN) and Parametric ReLU (PReLU) activation:
\begin{align*}
\mathbf{F}_k^{(1)}&={\rm PReLU}\left({\rm BN}(\mathbf{W}_{2k+1}*\mathbf{F}_k^{(0)}+\mathbf{b}_{2k+1})\right),\\
\mathbf{F}_k^{(0)}&=\mathbf{F}_k^{(0)}+{\rm BN}\left((\mathbf{W}_{2k+2}*\mathbf{F}_k^{(1)}+\mathbf{b}_{2k+2})\right),
\end{align*}
where $k=1,\dots,16$ indexes the residual blocks, and each $\mathbf{W}_{2k+1}$ and $\mathbf{W}_{2k+2}$ are $3 \times 3 \times m_2 \times m_2$ convolutional kernels. The residual connections within each block ensure that $\mathbf{F}_{k+1}^{(0)}$ retains information from previous layers, facilitating better gradient flow and improved learning stability.
After the residual blocks, we apply a convolutional layer with kernel size $3 \times 3$, followed by Batch Normalization and a large skip connection that adds $\mathbf{F}_1^{(0)}$ back to the result:
\begin{align*}
	\mathbf{F}_{\text{res}}=\mathbf{F}_1^{(0)}+{\rm BN}\left(\mathbf{W}_{\text{res}}*\mathbf{F}_{16}^{(0)}+\mathbf{b}_{\text{res}}\right),
\end{align*}
where $\mathbf{W}_{\text{res}}$ is a $3\times 3\times m_2\times m_2$ kernel.
Subsequently, the upsampling process begins. The output from the residual blocks, $\mathbf{F}_\text{res}$, is passed through two subpixel convolution blocks to achieve the desired spatial resolution. Each block doubles the spatial resolution while halving the number of channels using Pixel Shuffle operations, implemented as:
\begin{align*}
	\mathbf{F}_{\text{up}}^{(1)}&={\rm PixelShuffle}\left({\rm PReLU}(\mathbf{W}_{\text{up}}^{(1)}*\mathbf{F}_{\text{res}}+\mathbf{b}_{\text{up}}^{(1)}),r=2\right),\\
	\mathbf{F}_{\text{up}}^{(2)}&={\rm PixelShuffle}\left({\rm PReLU}(\mathbf{W}_{\text{up}}^{(2)}*\mathbf{F}_{\text{up}}^{(1)}+\mathbf{b}_{\text{up}}^{(2)}),r=2\right),
\end{align*}
where $\mathbf{W}_{\text{up}}^{(1)}$, $\mathbf{W}_{\text{up}}^{(2)}$ have dimensions $3 \times 3 \times 64 \times 64$, and $r$ denotes the upscale factor of 2 in each dimension. The Pixel Shuffle reorganizes the channels into spatial dimensions, thereby increasing the resolution of $\mathbf{F}_{\text{up}}^{(2)}$ to $m_2\times m_2$.
Finally, we apply a convolutional layer with kernel size $9\times9$, followed by a Tanh activation function to generate the final "high-resolution" output matrix:
\begin{align*}
	\mathbf{Y}_{\text{high}}={\rm Tanh}\left(\mathbf{W}_{\text{out}}*\mathbf{F}_{\text{up}}^{(2)}+\mathbf{b}_{\text{out}}\right),
\end{align*}
where $\mathbf{W}_{\text{out}}$ is a $9 \times 9 \times m_2 \times 1$ kernel. The Tanh activation ensures the output values are bounded, suitable for capturing the dispersion relation. The detailed structure of SRResNet as well as the workflow are shown in Figure \ref{SRResNet structure}.

Let $\{\omega_n^{m_1,j},\omega_n^{m_2,j})\}_{j=1,\cdots,\tilde{N}_{{\rm train}},n=1,\cdots,L}$ with $m_2=4\times m_1$ and $m_1\in\mathbb{N}$ denote the labeled training dataset, ${N}_{{\rm train}}:=L\tilde{N}_{{\rm train}}$ the total number of training data and
$\theta_2$ the parameters from SRResNet. To learn the mapping \eqref{F 2}, we use discrete mean square error as the loss function, 
\begin{align}\label{eq:loss-f2}
	\mathcal{L}_2(\theta_2)=\frac{1}{m^2\tilde{N}_{{\rm train}}}\sum_{j=1}^{\tilde{N}_{{\rm train}}}\sum_{n=1}^L\left\|G(\omega_n^{m_1,j};\theta_2)-\omega_n^{m_2,j}\right\|_F^2.
\end{align}
The dataset is also shuffled and the training is performed using batches of size $N_B$.

\section{Numerical tests}\label{Results}
We present in this section the performance of our proposed supervised learning approach. In these tests, we will use the mean relative error (MRE) as the primary evaluation metric to demonstrate the effectiveness of our model. MRE is particularly important in this context because it normalizes the error relative to the true values, providing a more intuitive understanding of the prediction accuracy in proportion to the magnitude of the true band function values. MRE is defined as follows:
\begin{align*}
{\rm MRE}:=\frac{1}{m^2N}\sum_{j=1}^{\tilde{N}_{\rm Test}}\sum_{n=1}^L\sum_{p,q=1}^m\frac{\left|\left(\hat{\omega}_n^{m,j}\right)_{pq}-\left(\omega_n^{m,j}\right)_{pq}\right|}{\left|\left(\omega_n^{m,j}\right)_{pq}\right|},
\end{align*}
where $\tilde{N}_{\rm Test}$ is the number of unit cells in the testing set, $L$ the number of band functions evaluated, $N:=L\tilde{N}_{\rm Test}$, $\hat{\omega}_n^{j}$ represents the predicted result from our model, $\omega_n^{j}$ is the true value derived from FEM.
\begin{figure}[hbt!]
	\centering
	\includegraphics[width = .8\textwidth,trim={0.8cm 0.8cm 1cm 0.5cm},clip]{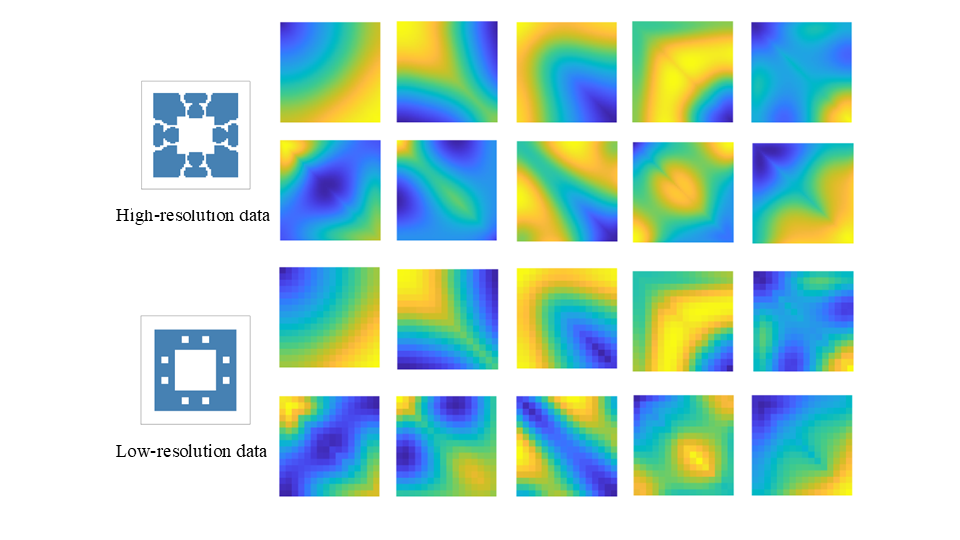}
	\centering
	\caption{Illustration of labeled data $\left((M_{\Omega,64}, I_{64}^{(n)}),\omega_{\text{H},n}^{64}\right)$ and $\left((M_{\Omega,16}, I_{16}^{(n)}),\omega_{\text{L},n}^{16}\right)$ for learning task \eqref{F 1} with $n=1,\cdots,10$: unit cell matrices $M_{\Omega,64}$ (top left), $M_{\Omega,16}$ (bottom left) and the first 10 band functions $\omega_{\text{H},n}^{64}$ (top right), $\omega_{\text{L},n}^{16}$ (bottom right), arranged from left to right and top to bottom. All figures are displayed using MATLAB built-in function \textit{imagesc}.}\label{data_set}
\end{figure}
\begin{figure}[hbt!]
	\centering
	\includegraphics[width = .8\textwidth,trim={0.8cm 0.8cm 1cm 0.5cm},clip]{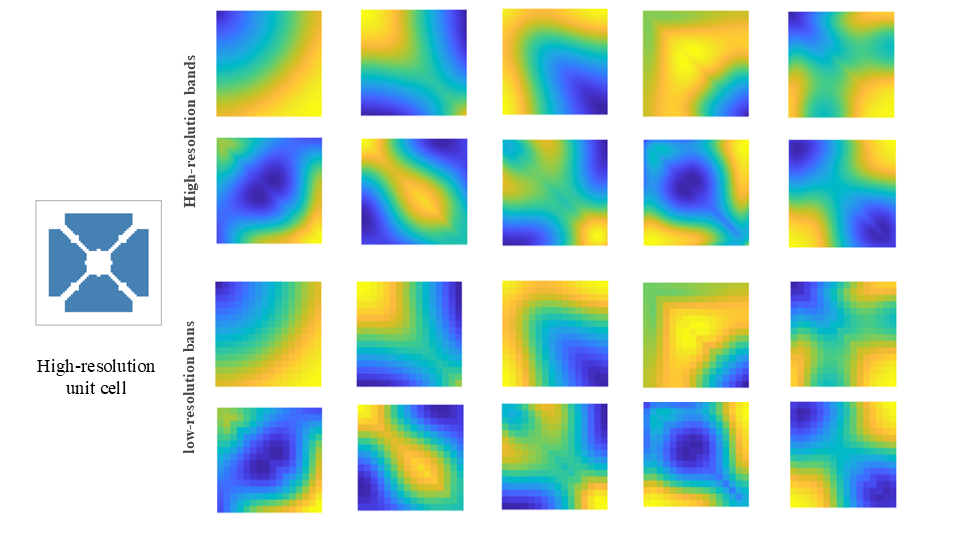}
	\centering
	\caption{Illustration of 10 labeled data $(\omega_{\text{H},n}^{16},\omega_{\text{H},n}^{64})$ for learning task \eqref{F 2} with $n=1,\cdots,10$: one given unit cell matrix $M_{\Omega,64}$ (left), the "high-resolution" band functions $\omega_{\text{H},n}^{64}$ (top right), and the "low-resolution" band functions $\omega_{\text{H},n}^{16}$ (bottom right), arranged from left to right and top to bottom. All figures are displayed using MATLAB built-in function \textit{imagesc}.}\label{data_set 2}
\end{figure}
\begin{figure}[hbt!]
	\centering
	\subfigure[Learning task \eqref{F 1} with "low-resolution" dataset without transfer learning]{\label{loss(1)}
		\includegraphics[width = .39\textwidth,trim={0.7cm 0.5cm 0.7cm 1.5cm},clip]{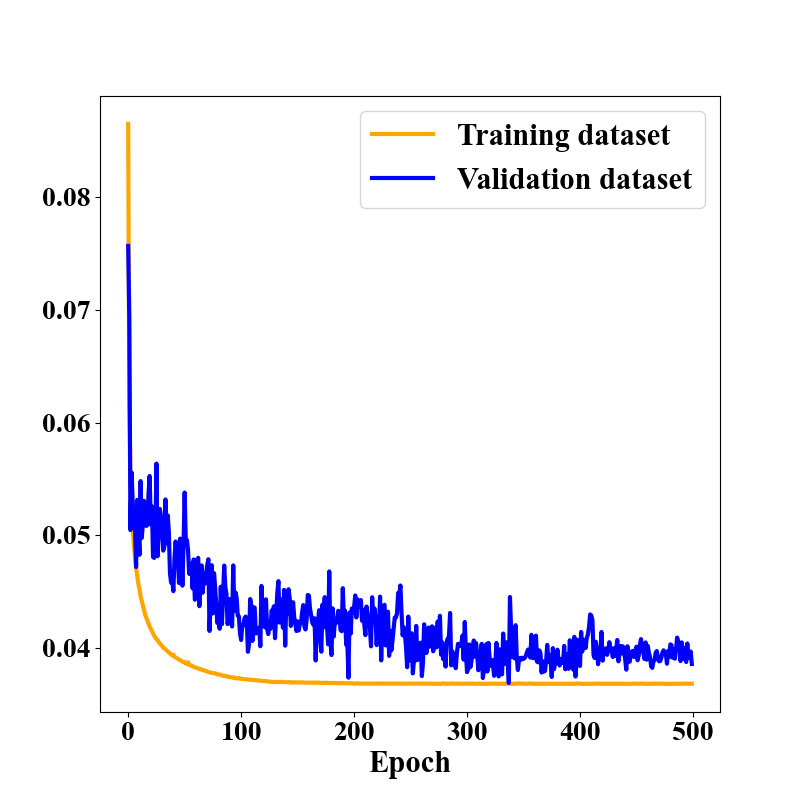}
	}%
	\subfigure[Learning task \eqref{F 1} with "high-resolution" dataset without transfer learning]{\label{loss(2)}
		\includegraphics[width = .39\textwidth,trim={0.7cm 0.5cm 0.7cm 1.5cm},clip]{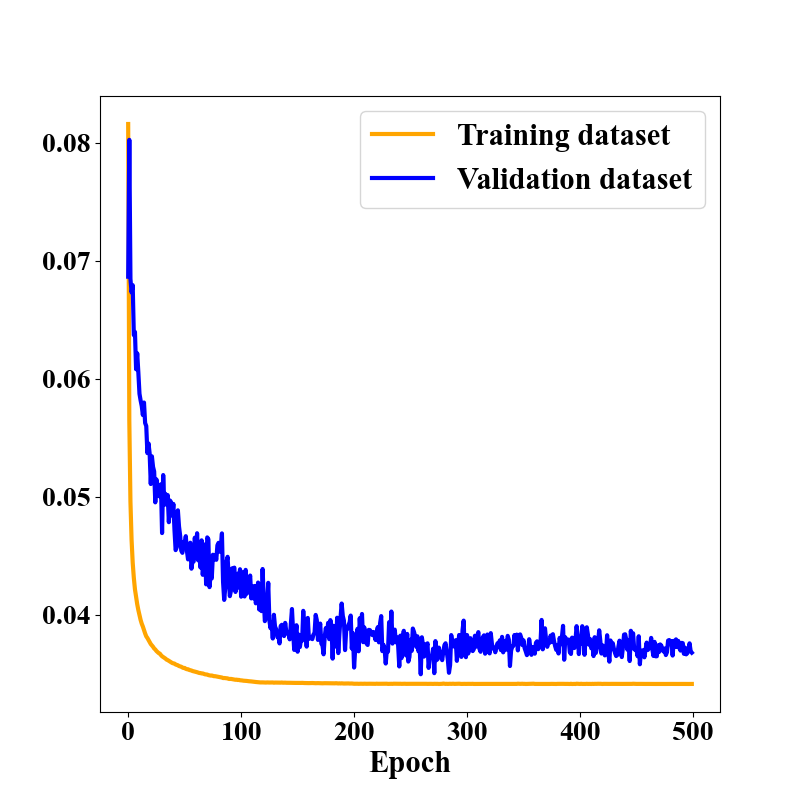}
	}%
\\
	\subfigure["low-resolution" bands 6-10 using randomly initialized model]{\label{loss(1) 1 5 _16}
		\includegraphics[width = .39\textwidth,trim={0.7cm 0.5cm 0.7cm 1.5cm},clip]{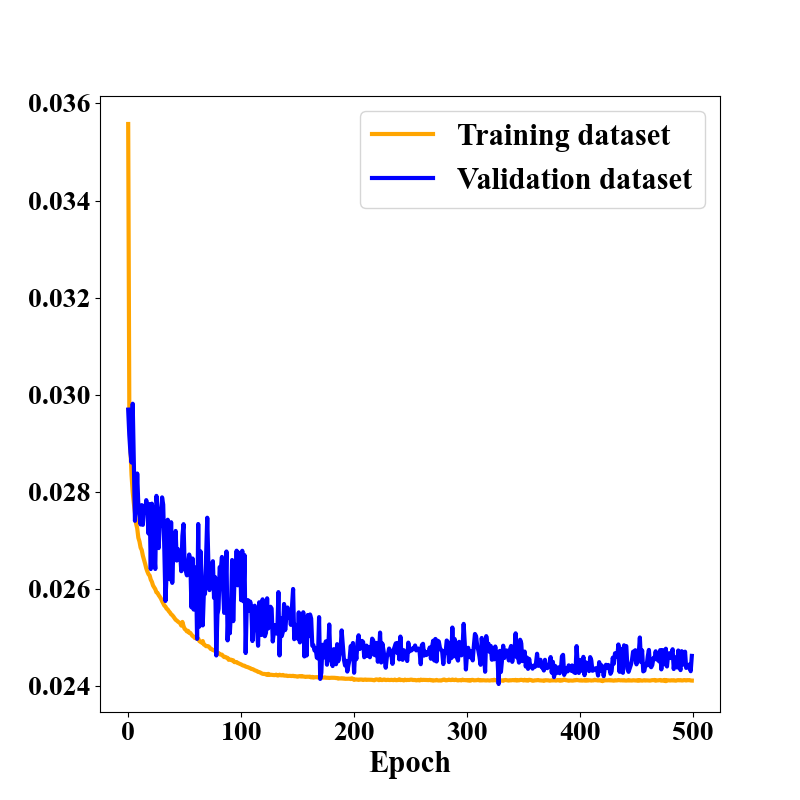}
	}%
	\subfigure["low-resolution" bands 6-10 using pre-trained model (trained on bands 1-5)]{\label{loss(2) 6 10 _16}
		\includegraphics[width = .39\textwidth,trim={0.7cm 0.5cm 0.7cm 1.5cm},clip]{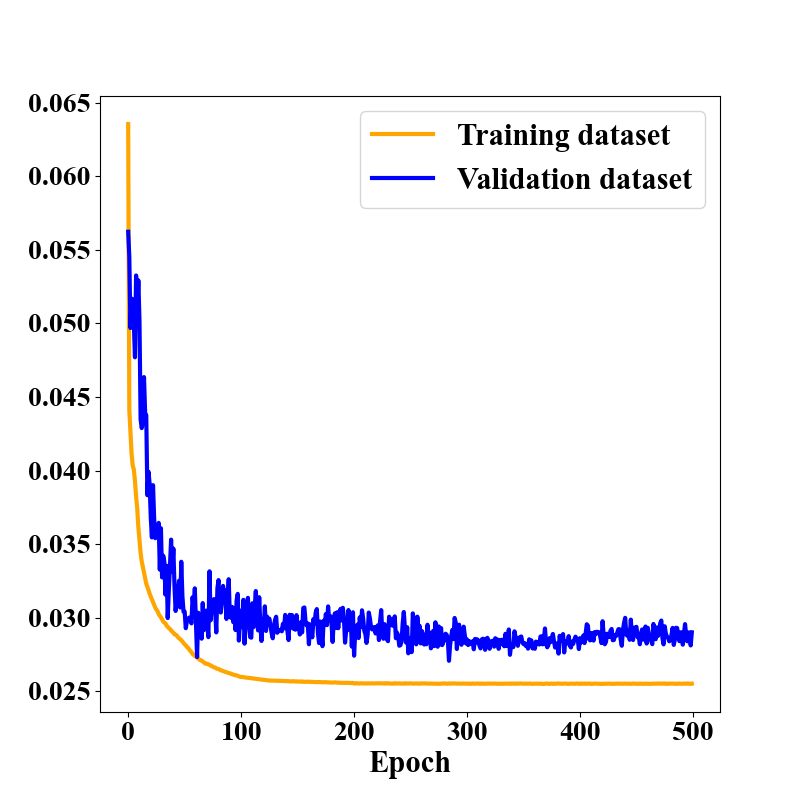}
	}%
\\
	\subfigure["high-resolution" bands 6-10 using randomly initialized model]{\label{loss(1) 1 5 _64}
		\includegraphics[width = .39\textwidth,trim={0.7cm 0.5cm 0.7cm 1.5cm},clip]{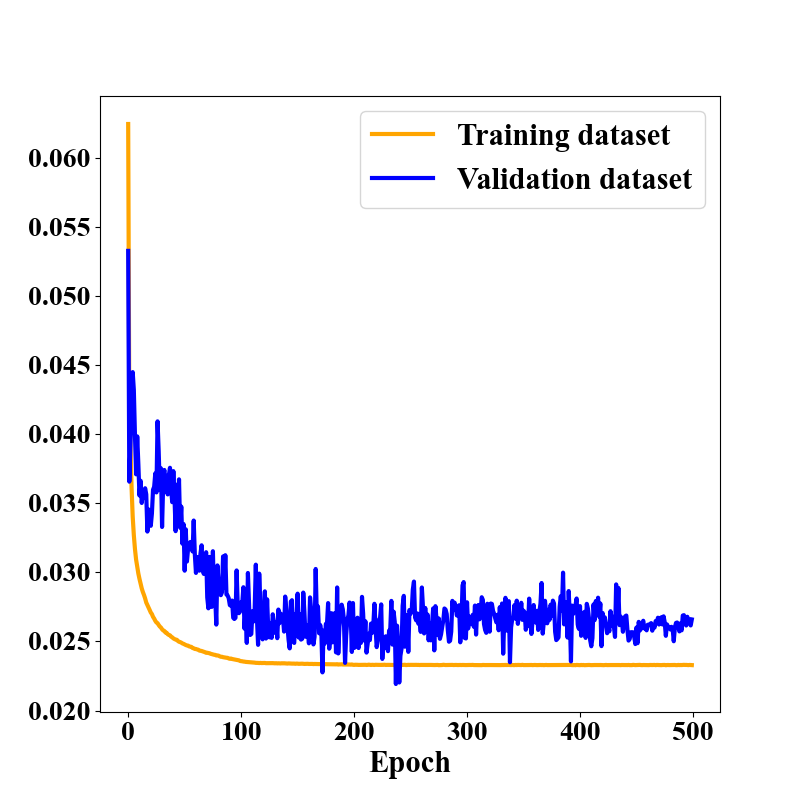}
	}%
	\subfigure["high-resolution" bands 6-10 using pre-trained model (trained on bands 1-5)]{\label{loss(2) 6 10 _64}
		\includegraphics[width = .39\textwidth,trim={0.7cm 0.5cm 0.7cm 1.5cm},clip]{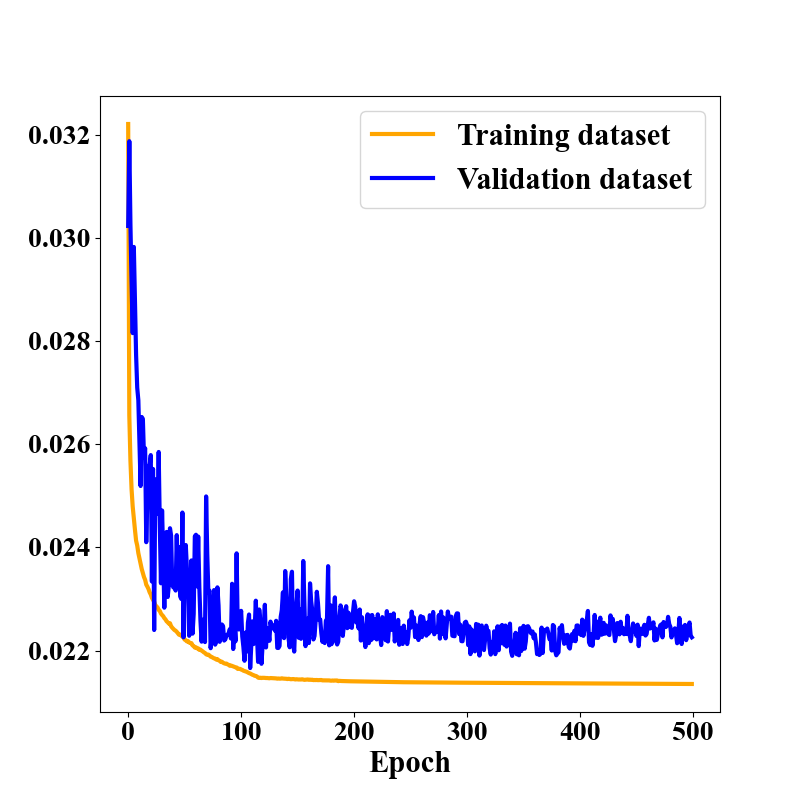}
	}%
	\centering
	\caption{Learning task \eqref{F 1}: MRE against training epochs.} 
\label{loss 3}
\end{figure}
\begin{table}[hbt!]
	\centering
	\begin{tabular}{|c|c|c|c|c|c|c|}
		\hline
		\multicolumn{6}{|c|}{\textbf{Bands 1-10}} & \textbf{MRE} \\ \hline
		Band number & 1 & 2 & 3 & 4 & 5  & \multirow{4}{*}{4.68\%}\\ \cline{1-6}
		MRE & 4.78\%&4.67\%&4.51\%&4.53\%&4.62\%& \\ \cline{1-6}
		Band number & 6 & 7 & 8 & 9 & 10& \\ \cline{1-6}
		MRE & 4.41\%&4.52\%&4.87\%&4.91\%&4.98\%& \\ \hline
		\multicolumn{6}{|c|}{\textbf{Bands 6-10, randomly initialized model}}  & \textbf{MRE} \\ \hline
		Band number & 6 & 7 & 8 & 9 & 10& \multirow{2}{*}{3.62\%} \\ \cline{1-6}
		MRE & 3.68\% & 3.67\% & 3.58\% & 3.52\% & 3.64\%& \\ \hline
		\multicolumn{6}{|c|}{\textbf{Bands 6-10, pre-trained model}} & \textbf{MRE}\\ \hline
		Band number & 6 & 7 & 8 & 9 & 10 & \multirow{2}{*}{2.86\%} \\ \cline{1-6}
		MRE & 2.74\%&2.87\%&2.78\%&2.92\%&2.98\%& \\ \hline
	\end{tabular}
	\caption{Outcome of learning task \eqref{F 1} with "low-resolution" dataset.}
	\label{tab1}
\end{table}
\begin{table}[hbt!]
	\centering
	\begin{tabular}{|c|c|c|c|c|c|c|}
		\hline
		\multicolumn{6}{|c|}{\textbf{Bands 1-10}}& \textbf{MRE} \\ \hline
		Band number & 1 & 2 & 3 & 4 & 5 & \multirow{4}{*}{$4.78\%$}\\ \cline{1-6}
		MRE & $4.92\%$&$4.77\%$&$4.68\%$&$4.61\%$&$4.76\%$& \\ \cline{1-6}
		Band number & 6 & 7 & 8 & 9 & 10& \\ \cline{1-6}
		MRE & $4.59\%$&$4.74\%$&$4.78\%$&$4.93\%$&$4.98\%$& \\ \hline
		\multicolumn{6}{|c|}{\textbf{Bands 6-10, randomly initialized model}}& \textbf{MRE} \\ \hline
		Band number & 6 & 7 & 8 & 9 & 10 & \multirow{2}{*}{3.59\%}\\ \cline{1-6}
		MRE & 3.57\%&3.62\%&3.61\%&3.53\%&3.61\%& \\ \hline
		\multicolumn{6}{|c|}{\textbf{Bands 6-10, pre-trained model}}& \textbf{MRE} \\ \hline
		Band number & 6 & 7 & 8 & 9 & 10 & \multirow{2}{*}{2.75\%}\\ \cline{1-6}
		MRE & 2.51\%&2.68\%&2.74\%&2.91\%&2.92\%& \\ \hline
	\end{tabular}
	\caption{Outcome of learning task \eqref{F 1} with "high-resolution" dataset.}
	\label{tab2}
\end{table}

\begin{figure}[hbt!]
	\centering
	\subfigure["low-resolution" bands 1-10]{\label{tong ji 1}
		\includegraphics[width = .4\textwidth,trim={0cm 0cm 0cm 0cm},clip]{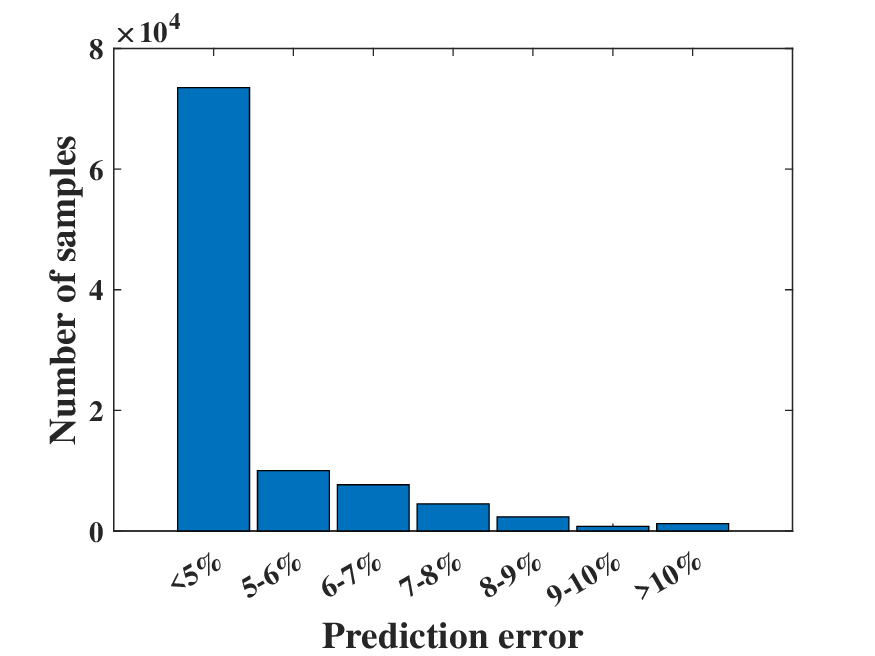}
	}%
	\subfigure["high-resolution" bands 1-10]{\label{tong ji 2}
		\includegraphics[width = .4\textwidth,trim={0cm 0cm 0cm 0cm},clip]{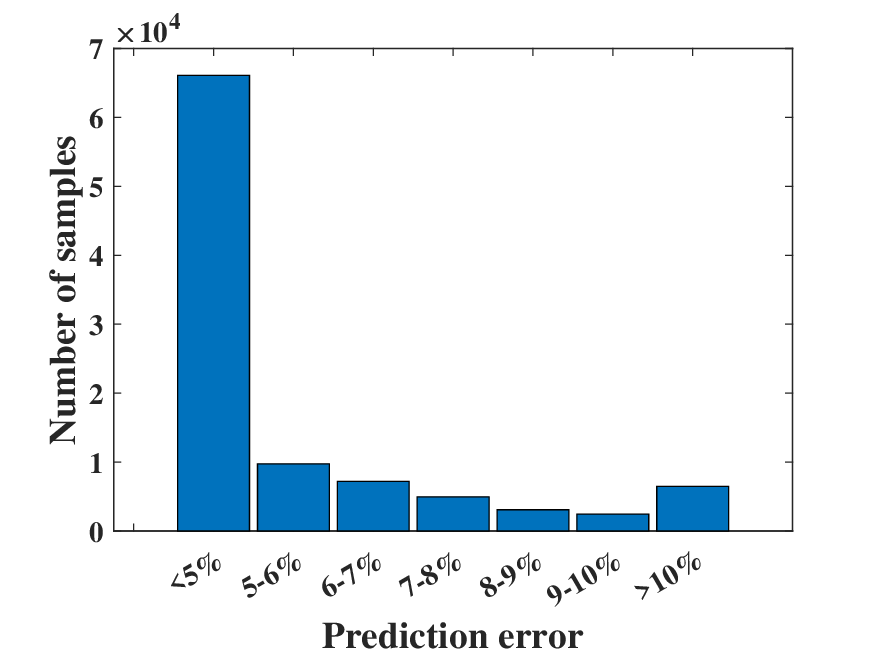}
	}%
	\\
	\subfigure["low-resolution" bands 6-10 using randomly initialized model]{\label{tong ji 5}
		\includegraphics[width = .4\textwidth,trim={0cm 0cm 0cm 0cm},clip]{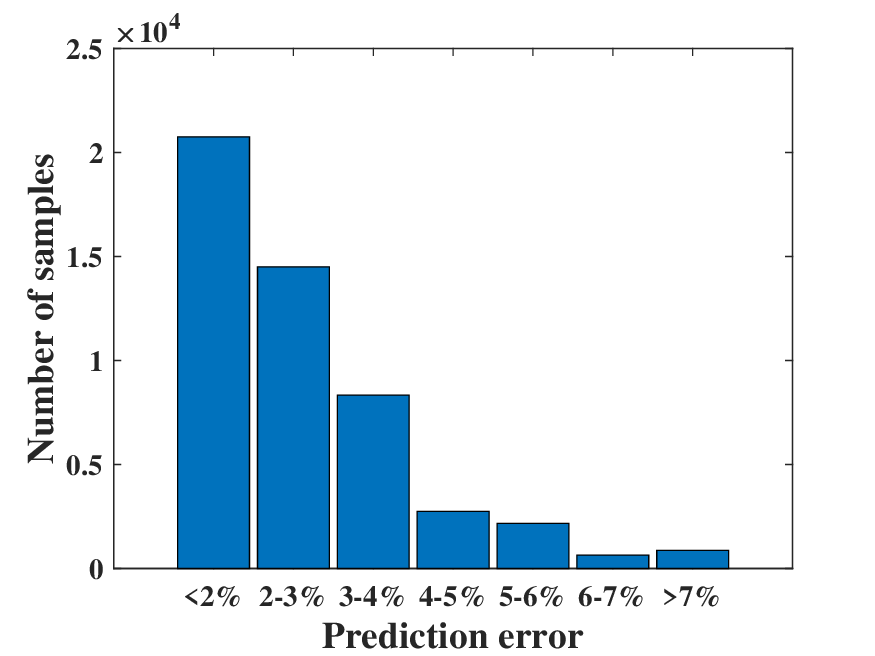}
	}%
	\subfigure["low-resolution" bands 6-10 using pre-trained model (trained on bands 1-5)]{\label{tong ji 6}
		\includegraphics[width = .4\textwidth,trim={0cm 0cm 0cm 0cm},clip]{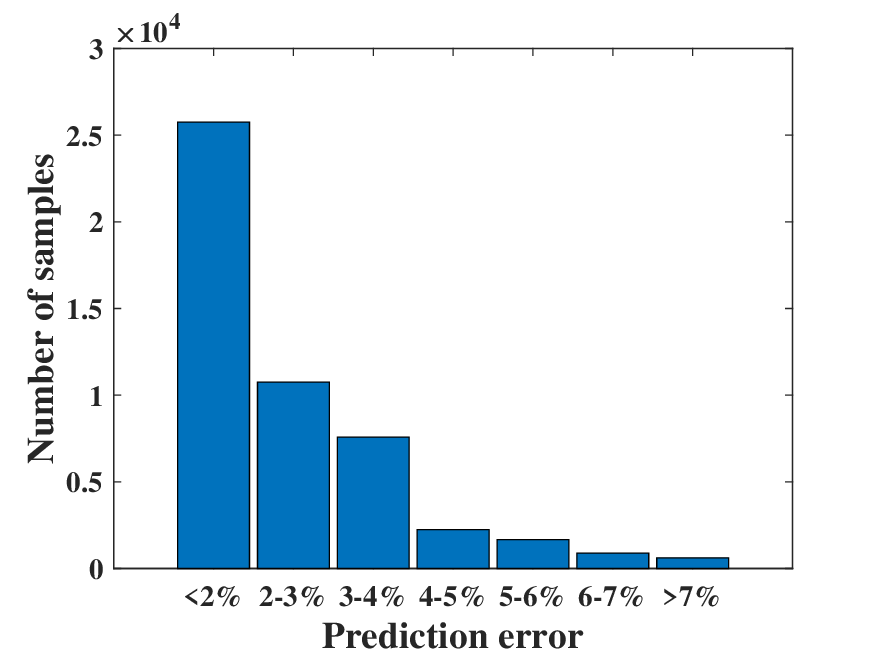}
	}%
	\\
	\subfigure["high-resolution" bands 6-10 using randomly initialized model]{\label{tong ji 3}
		\includegraphics[width = .4\textwidth,trim={0cm 0cm 0cm 0cm},clip]{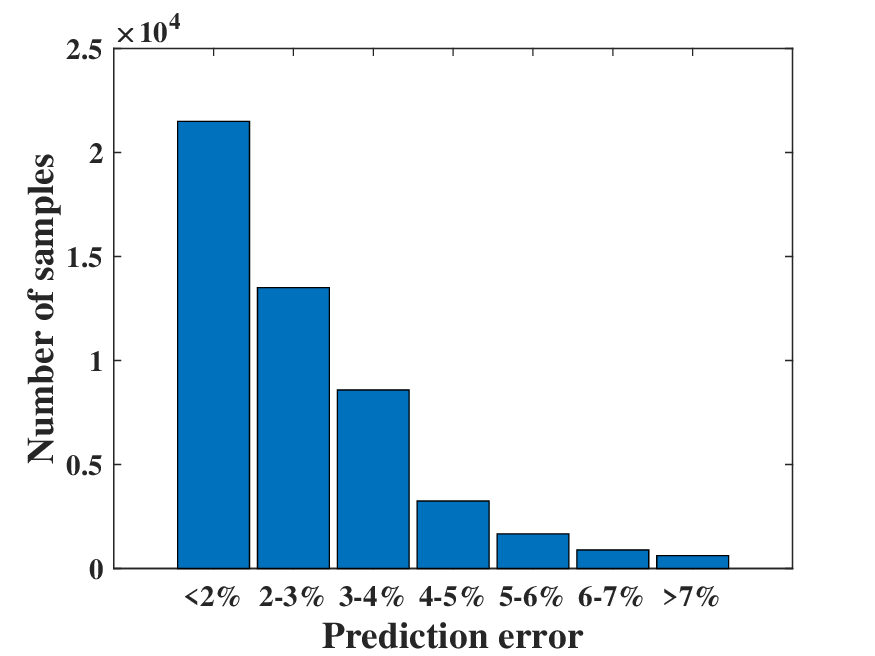}
	}%
	\subfigure["high-resolution" bands 6-10 using pre-trained model (trained on bands 1-5)]{\label{tong ji 4}
		\includegraphics[width = .4\textwidth,trim={0cm 0cm 0cm 0cm},clip]{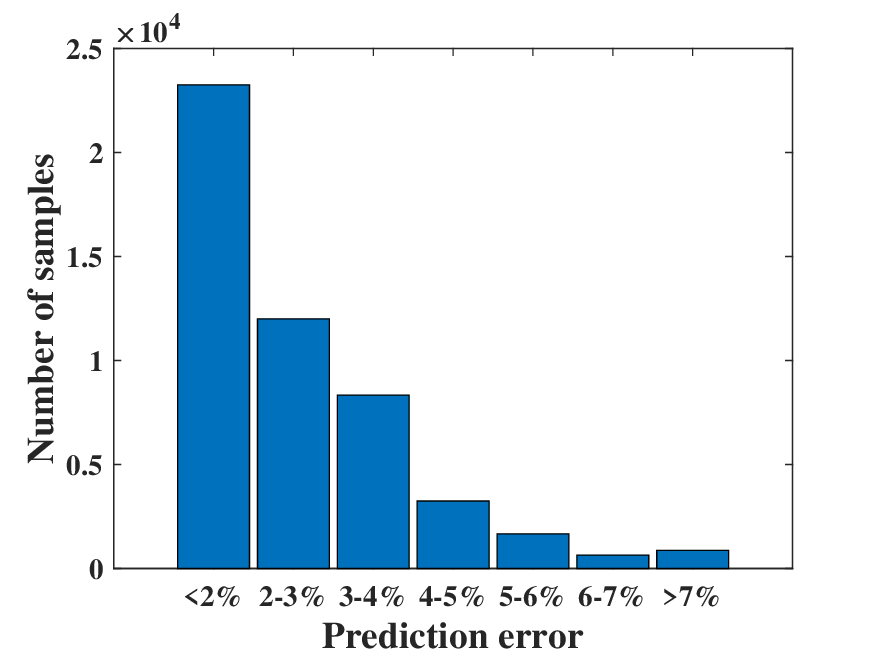}
	}%
	\centering
	\caption{U-net prediction distribution.}
\label{tong ji}
\end{figure}
\subsection{Training data}\label{Data preparation}
First, we introduce the generation of training data for our proposed supervised learning scheme. 



For each unit cell $\Omega$ and its discretization $M_{\Omega,m}$ \eqref{eq:dis-ucell}, we obtain the labeled data $((M_{\Omega,m},I_m^{n}),\omega_n^{m}(\mathbf{k}))$ for $n=1,\cdots, L$ by solving the Helmholtz eigenvalue problem \eqref{simply}, parameterized by the wave vector $\mathbf{k} \in \mathcal{B}$. For a given wave vector, we employ the conforming Galerkin Finite Element Method as outlined in \cite{wang2023dispersion,wang2023hp}. Each unit cell is discretized using triangular elements, with each pixel divided into two isosceles right triangles. For "high-resolution" unit cells with $m=64$, the FEM solution results in a degree of freedom (DoF) of $65^2$. In contrast, "low-resolution" unit cells with $m=16$ have a DoF of $17^2$. To maintain consistency in input-output dimensions within the U-Net, we calculate the dispersion relations for each $64 \times 64$ unit cell over a $64 \times 64$ grid of uniformly distributed wave vectors in the first Brillouin zone, encompassing the first 10 eigenfrequencies, i.e., $L=10$. These "high-resolution" dispersion relations are represented as $\omega^{64}_n$ for $n=1,\cdots, 10$. Additionally, we compute the "low-resolution" dispersion relations over a $16 \times 16$ grid of uniformly distributed wave vectors, resulting in $\omega^{16}_{\text{H},n}$ for $n=1,\cdots, 10$.
We also compute another set of "low-resolution" dispersion matrices for each $M_{\Omega,16}$ over $16 \times 16$ wave vectors, which are denoted as $\omega^{16}_{\text{L},n}$ for $n=1,\cdots, 10$.

This process yields three distinct datasets. For training U-net (Section \ref{UNET}), to achieve the mapping \eqref{F 1}, we utilize two datasets: the "high-resolution" dataset consists of $1,000,000$ pairs of $\left((M_{\Omega_j,64}, I_{64}^{(n)}),\omega_{\text{H},n}^{64,j}\right)$, for $j=1,\cdots,100,000$, $n=1,\cdots,10$. Meanwhile, the "low-resolution" dataset contains $1,000,000$ pairs of $\left((M_{\Omega_j,16}, I_{16}^{(n)}),\omega_{\text{L},n}^{16,j}\right)$, for $j=1,\cdots,100,000$, $n=1,\cdots,10$. Figure \ref{data_set} illustrates a sample visualization of the dataset. For training SRResNet (Section \ref{Super-Resolution Residual Network}) to achieve mapping \eqref{F 2}, we have a separate dataset consisting of pairs $\left(\omega_{\text{H},n}^{16,j},\omega_{\text{H},n}^{64,j}\right)$, for $j=1,\cdots,100,000$, $n=1,\cdots,10$, serving as the ground truth for training and evaluating our deep learning model. Figure \ref{data_set 2} illustrates a sample visualization of this dataset. We split this dataset into training, validation, and testing sets, comprising $80\%$, $10\%$, and $10\%$ of the samples, respectively. During the training process, we shuffle the training and testing sets and choose batch size $N_B=16$.
\begin{figure}[hbt!]
	\centering
	\includegraphics[width = .85\textwidth,trim={0.8cm 0.7cm 1cm 0.5cm},clip]{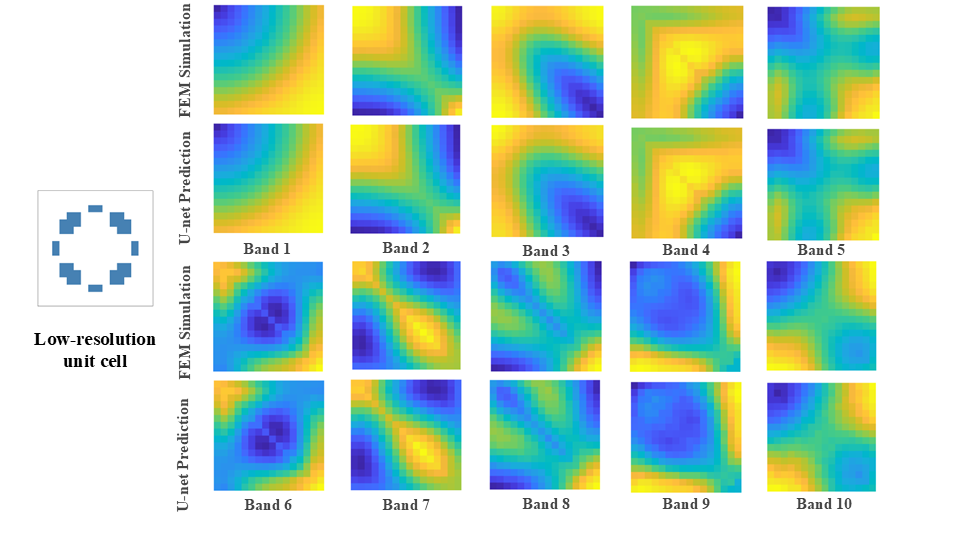}
	\centering
	\caption{Performance of the learned model from learning task \eqref{F 1} using "low-resolution" dataset.}\label{p 3}
\end{figure}
\begin{figure}[h]
	\centering
	\includegraphics[width = .85\textwidth,trim={0.8cm 0.5cm 1cm 0.5cm},clip]{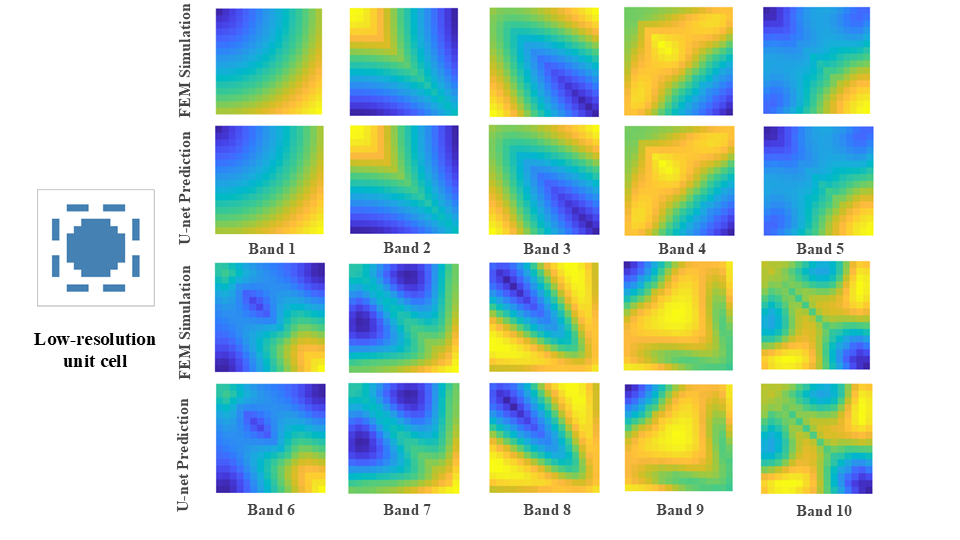}
	\centering
	\caption{Performance of the learned model from learning task \eqref{F 1} using "low-resolution" dataset with transfer learning (bands 1-5 using randomly initialized model and bands 6-10 using pre-trained model).}\label{p 5}
\end{figure}
\subsection{Learning task \eqref{F 1}}

To train the parameters by minimizing the loss function \eqref{eq:loss-f1}, we employ the Stochastic Gradient Descent (SGD) method with an initial learning rate of $0.01$. To ensure stable convergence, the learning rate is reduced by a factor of 10 after every 100 epochs. Once well trained, the U-Net can rapidly predict the dispersion relation, within an acceptable accuracy range, of the input unit cell in batches several orders of magnitude faster than the traditional numerical simulations. 
\begin{figure}[hbt!]
	\centering
	\includegraphics[width = .85\textwidth,trim={0.8cm 0.5cm 1cm 0.5cm},clip]{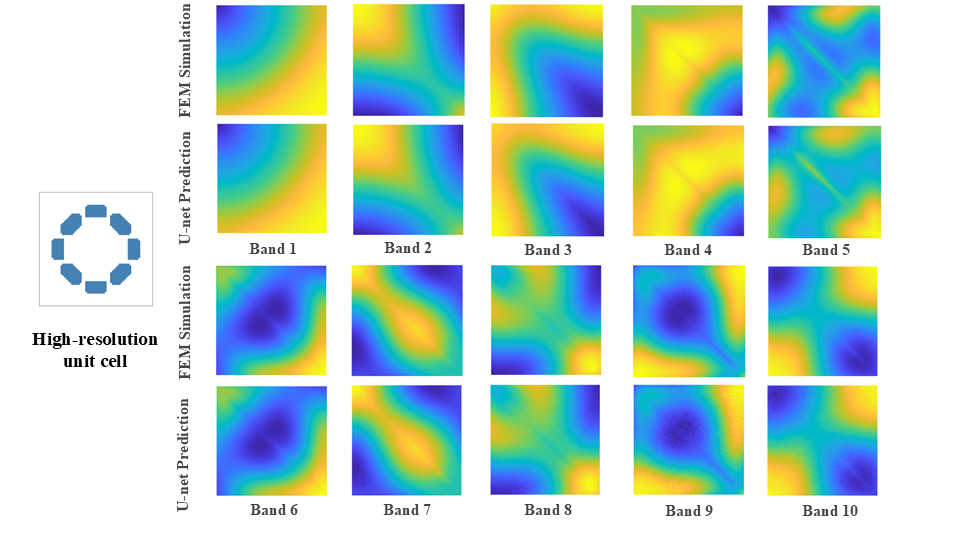}
	\centering
	\caption{Performance of the learned model from learning task \eqref{F 1} using "high-resolution" dataset.}\label{p 2}
\end{figure}

\begin{figure}[hbt!]
	\centering
	\includegraphics[width = .85\textwidth,trim={0.8cm 0.7cm 1cm 0.5cm},clip]{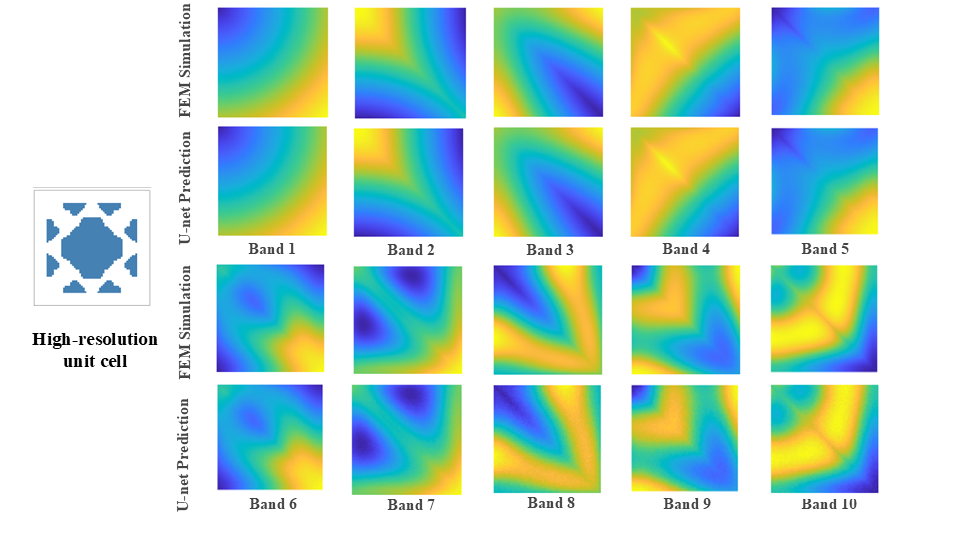}
	\centering
	\caption{Performance of the learned model from learning task \eqref{F 1} using "high-resolution" dataset with transfer learning (bands 1-5 using randomly initialized model and bands 6-10 using pre-trained model).}\label{p 4}
\end{figure}
Figure \ref{loss 3} illustrates the progression of the mean relative error across training and validation datasets as the number of epochs increases. The distribution of prediction errors is shown in Figure \ref{tong ji}. Initially, we evaluated the U-Net model’s capability to jointly predict the first ten band functions. Figures \ref{loss(1)} and \ref{loss(2)} display the results for both “low-resolution” and “high-resolution” unit cells. It is evident that the model converges well over time, achieving mean relative errors below $5\%$ for both resolutions, as summarized in Tables \ref{tab1} and \ref{tab2}.

\begin{figure}[hbt!]
	\centering
	\subfigure[MRE against the training epochs]{\label{loss(1) super}
		\includegraphics[width = .34\textwidth,trim={0.3cm 0.7cm 2.7cm 3cm},clip]{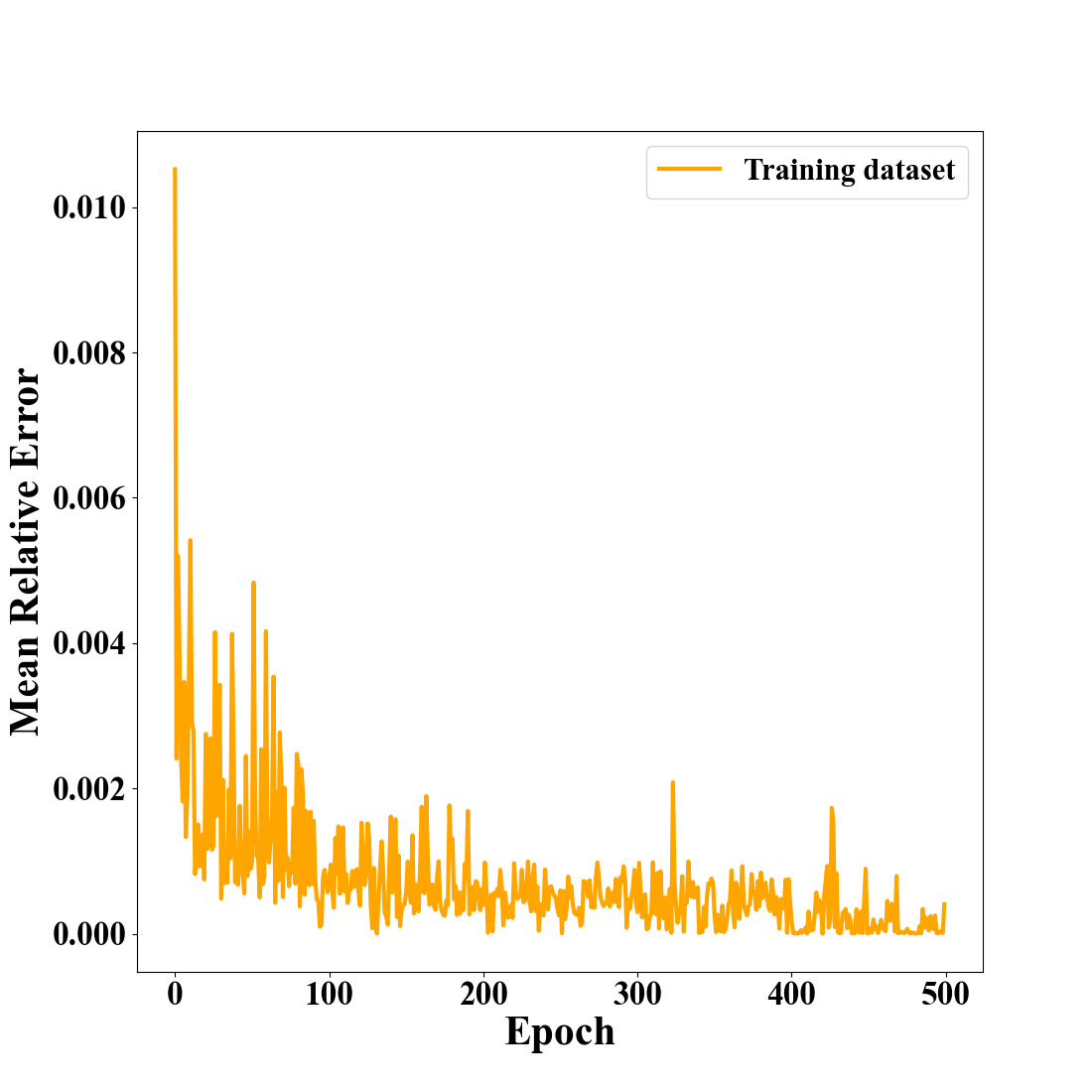}
	}%
	\subfigure[Prediction distribution]{\label{super tong ji}
		\includegraphics[width = .45\textwidth,trim={0cm 0cm 0cm 0cm},clip]{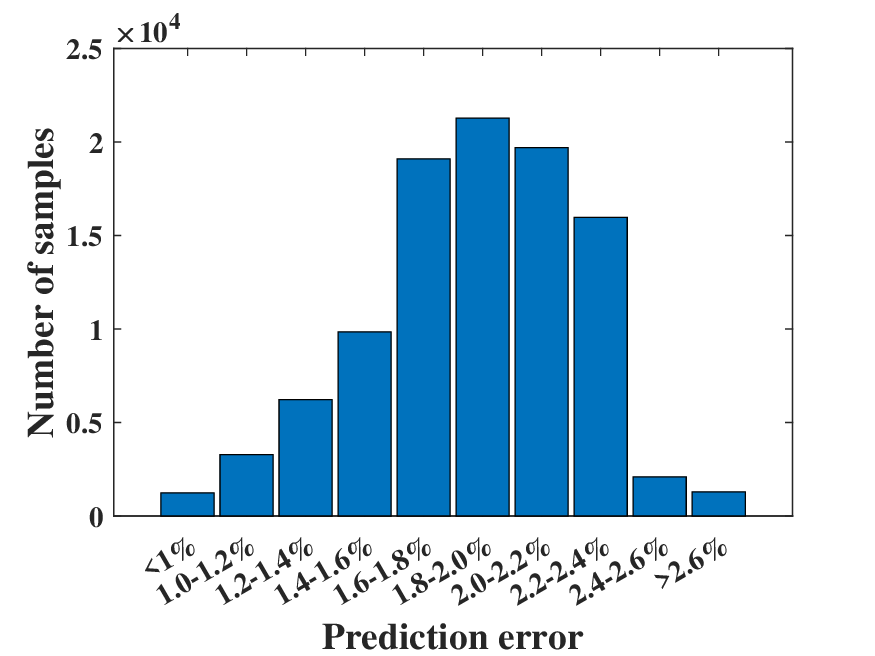}
	}
	\centering
	\caption{Learning task \eqref{F 2}.}
\label{loss 3 super}
\end{figure}

\begin{table}[hbt!]
	\centering
	\begin{tabular}{|c|c|c|c|c|c|c|}
		\hline
		\multicolumn{6}{|c|}{\textbf{Bands 1-10}}& \textbf{MRE} \\ \hline
		Band number & 1 & 2 & 3 & 4 & 5 & \multirow{4}{*}{$1.85\%$}\\ \cline{1-6}
		MRE & $1.84\%$&$1.79\%$&$1.77\%$&$1.83\%$&$1.88\%$& \\ \cline{1-6}
		Band number & 6 & 7 & 8 & 9 & 10& \\ \cline{1-6}
		MRE & $1.82\%$&$1.89\%$&$1.89\%$&$1.90\%$&$ 1.87\%$& \\ \hline
	\end{tabular}
	\caption{Outcome of learning task \eqref{F 2}.}
	\label{tab3}
\end{table}

\begin{table}[hbt!]
	\centering
	\begin{tabular}{|c|c|c|c|c|c|c|}
		\hline
		\multicolumn{6}{|c|}{\textbf{Bands 1-10}}& \textbf{MRE} \\ \hline
		Band number & 1 & 2 & 3 & 4 & 5 & \multirow{4}{*}{$6.88\%$}\\ \cline{1-6}
		MRE & $ 5.94\%$&$ 5.91\%$&$ 6.23\%$&$ 5.93\%$&$ 6.46\%$& \\ \cline{1-6}
		Band number & 6 & 7 & 8 & 9 & 10& \\ \cline{1-6}
		MRE & $ 8.31\%$&$ 7.35\%$&$ 8.01\%$&$ 7.42\%$&$  7.21\%$& \\ \hline
	\end{tabular}
	\caption{Linear interpolation}
	\label{linear interpolation}
\end{table}

We further examine the U-Net’s performance under a transfer learning scenario. Figures \ref{loss(1) 1 5 _16} and \ref{loss(1) 1 5 _64} present the training results for bands 6 to 10 on both "low-resolution" and "high-resolution" unit cells, using a randomly initialized U-Net model. As outlined in Tables \ref{tab1} and \ref{tab2}, the mean relative errors are $3.62\%$ and $3.59\%$, respectively. Subsequently, we re-trained the U-Net model initialized on bands 1 to 5 as a pre-trained model, employing only half the training and validation samples. The results, shown in Figures \ref{loss(2) 6 10 _16} and \ref{loss(2) 6 10 _64}, indicate improved predictive accuracy with mean relative errors dropping to $2.86\%$ and $2.75\%$. This highlights the advantage of transfer learning, allowing efficient training for higher-index bands without needing to fully retrain the model from scratch. 

Figures \ref{p 3}-\ref{p 4} show random examples comparing U-Net predictions to FEM simulations for bands 1 through 10. The results, visualized with MATLAB's \textit{imagesc} function, demonstrate good alignment between predicted dispersion relations and ground truth data. While prediction accuracy is lower for the first and last bands, likely due to their distinct structural characteristics and broader frequency ranges, the U-Net model performs robustly overall. The model successfully captures key properties of PhCs, including bandgap information. These results suggest that the U-Net can reliably predict PhCs’ dispersion relations, supporting further analysis of properties like phase and group velocities.

\begin{figure}[hbt!]
	\centering
	\subfigure{\label{super pre1}
		\includegraphics[width = .82\textwidth,trim={0.8cm 2cm 1cm 1cm},clip]{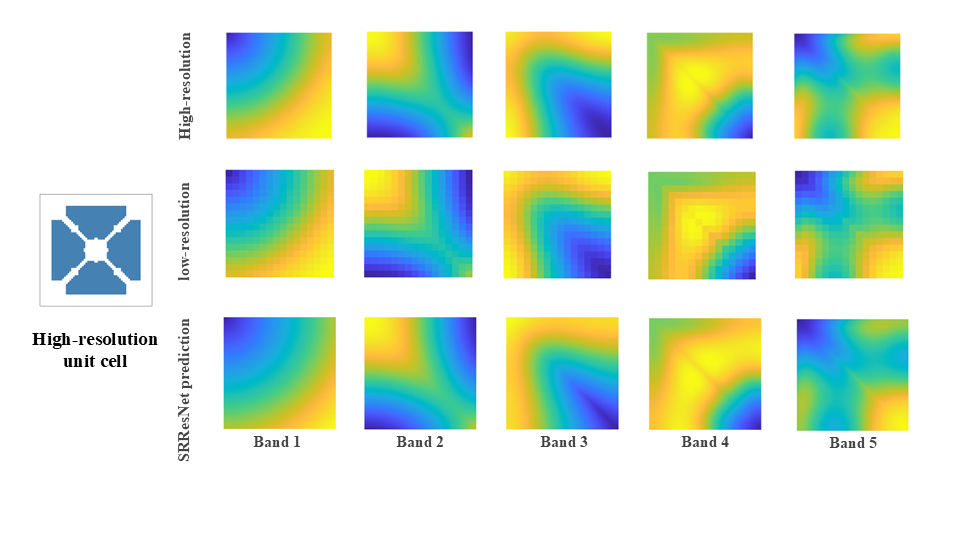}
	}%
	\\
	\subfigure{\label{super pre2}
		\includegraphics[width = .85\textwidth,trim={0.8cm 2.2cm 1cm 1cm},clip]{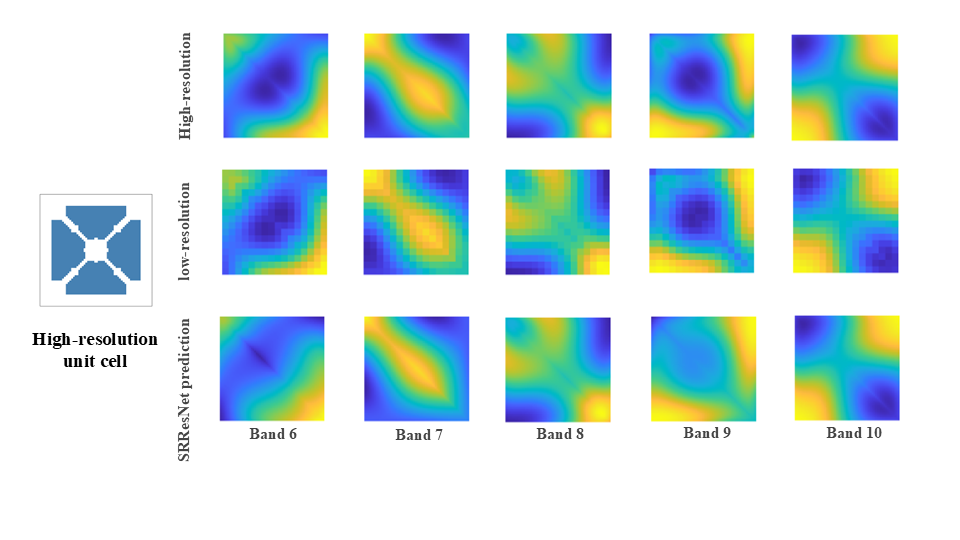}
	}%
	\centering
	\caption{Performance of the learned model from learning task \eqref{F 2}.}\label{super pre}
\end{figure}
\subsection{Learning task \eqref{F 2}}

To train SRResnet for the supervised learning task \eqref{F 2} by minimizing the loss function \eqref{eq:loss-f2}, we use Stochastic Gradient Descent (SGD) with an initial learning rate of $0.01$. The learning rate is reduced by a factor of 10 every 100 epochs, which enhances convergence and improves the accuracy of the super-resolved dispersion relations.

Figure \ref{loss(1) super} depicts the mean relative errors throughout the training process, demonstrating how the model converges over increasing epochs. The distribution of prediction errors is shown in Figure \ref{super tong ji}, further validating the model’s performance. As summarized in Table \ref{tab3}, the mean relative error between the predicted band functions and the simulated dispersion relation remains consistently below $2\%$.

A representative comparison between the dispersion relations predicted by SRResNet and those obtained from FEM simulations is presented in Figure \ref{super pre}. We focused on predicting the top ten band functions, which were visualized using MATLAB's \textit{imagesc} function. The close alignment between the predicted dispersion relations and the FEM ground truth confirms that the SRResNet effectively enhances the resolution of the band functions. This result implies that we can achieve detailed insights into a given band function using fewer parameter points within the domain, significantly improving efficiency.

Furthermore, we compare performance of the learned model from learning task \eqref{F 2} with a linear interpolation approach, which can also generate fine scale information from coarse scale. As shown in Table \ref{linear interpolation}, the mean relative error for linear interpolation is $6.88\%$, considerably higher than the error achieved by the learned model as reported in Table \ref{tab3}. This highlights the superior accuracy and efficacy of our proposed method over traditional interpolation techniques for resolution enhancement in photonic crystal band function prediction.

\section{Conclusion}\label{sec:conclusion}




In conclusion, this work demonstrates the successful application of deep learning models, specifically U-Net combined with transfer learning and Super-Resolution techniques, for predicting the dispersion relations of 2D photonic crystals across the Brillouin zone. By overcoming the limitations of traditional numerical methods , such as high computational costs and the requirement of fine mesh generation, our proposed model efficiently predicts high-resolution band structures from low-resolution data obtained through the Finite Element Method. Furthermore, our model's ability to predict the entire dispersion relation in a unified framework, rather than treating each band function separately, offers a more cohesive and reliable tool for the analysis of photonic crystal structures. 

The results confirm that our approach provides accurate predictions for the first several bands functions of 2D PhCs, with notable improvements in computational efficiency compared to conventional methods. This innovation not only accelerates the design and optimization process for photonic crystals but also broadens the potential applications of machine learning in complex physical systems. Overall, the synergy between advanced machine learning models and well-established numerical methods highlights a promising direction for future research and development in photonic crystal design, potentially benefiting the broader field of computational physics and materials science.

\printbibliography
\end{document}